\documentclass[12pt]{amsart}

\usepackage{amsmath,amssymb,enumerate}
\usepackage{relsize}%

\usepackage{epsfig,fancyhdr,color}

\usepackage{amssymb}
\usepackage{amsmath,amsthm}
\usepackage{bbm}
\usepackage{latexsym}
\usepackage{amscd}
\usepackage{psfrag}
\usepackage{graphicx}
\usepackage[all]{xy}

\usepackage{rotating}
\usepackage{adjustbox}

\usepackage{xcolor}

\usepackage[font=small,labelfont=bf]{caption}

\usepackage{mathtools}

\usepackage{tikz-cd}

\usepackage{stackrel}

\usepackage[utf8]{inputenc}
\usepackage{amsmath}
\usepackage{amssymb}

\usepackage{mathabx,epsfig}
\def\acts{\mathrel{\reflectbox{$\righttoleftarrow$}}}

\usepackage[utf8]{inputenc}

\numberwithin{equation}{section}

\definecolor{NoteColor}{rgb}{1,0,0}

\newtheorem{theorem}{\rm\bf Theorem}[section]

\newtheorem{lemma}[theorem]{\rm\bf Lemma}
\newtheorem{corollary}[theorem]{\rm\bf Corollary}
\newtheorem*{theorem 1}{\rm\bf Proposition 1}
\newtheorem*{theorem 2}{\rm\bf Proposition 2}

\theoremstyle{definition}

\theoremstyle{remark}
\newtheorem{remark}[theorem]{\rm\bf Remark}

\def\interieur#1{\mathord{\mathop{\kern 0pt #1}\limits^\circ}}

\def\interieur#1{\mathord{\mathop{\kern 0pt #1}\limits^\circ}}

\def\hyperp{{\rm I}\kern-.3ex{\rm H}}

\def\hyperl{\mathbb H}

\usepackage{mathabx,epsfig}
\def\acts{\mathrel{\reflectbox{$\righttoleftarrow$}}}

\newcommand{\OSp}{\operatorname{OSp}}
\newcommand{\osp}{\operatorname{OSp}(1|2)}

\newcommand{\Vol}{\operatorname{Vol}}

\begin{document}

\title[3d Super Hyperbolic Geometry]{3d Super Hyperbolic Geometry}
\
\author{Robert Penner} \address {\hskip -2.5ex Institut des Hautes \'Etudes Scientifiques\\ 35 route des Chartres\\ Le Bois Marie\\ 91440 Bures-sur-Yvette\\ France\\ {\rm and}~Mathematics Department, UCLA\\ Los Angeles, CA 90095\\USA} \email{rpenner{\char'100}ihes.fr}


 \date{\today}


\begin{abstract}
Wigner's unitary representation of the Lorentz group
is extended to a representation of 
the complex orthosymplectic Lie super group ${\rm OSp}_{\mathbb C}(1{\rm |}2)$ acting on
Minkowski (3,1${\rm |}$4)-dimensional super space essentially by Hermitean conjugation. The invariant quadratic form is
$x_1x_2-x\bar x+\phi\psi+\bar\phi\bar\psi$ in Wigner's coordinates
$x_1,x_2\in{\mathbb R}$ and $x\in{\mathbb C}$, 
where $\phi,\psi$ are Dirac fermions.  The extended action
is linear in the super space variables, but not quadratic in 
the odd group variables, and is described by an explicit even
purely imaginary auxiliary parameter defined on the product of the Lie super group 
with its Lie algebra.  This extension of Wigner's representation opens the door 
to studying geometry in super hyperbolic three-space, initiated here with foundations on super geodesics, triangles and  tetrahedra, and culminating in a proof of divergence of the volume of a typical ideal tetrahedron and a discussion of the non-zero fermionic correction to the Schl\"afli formula.
\end{abstract}

\maketitle





{\let\thefootnote\relax\footnote{{{Keywords: super three-dimensional hyperbolic geometry, complex ${\mathcal N}=1$ orthosymplectic group, Wigner's unitary representation of the Lorentz group, super three-dimensional Schl\"afli formula} }}}

{\let\thefootnote\relax\footnote{{{It is a pleasure to thank Igor Frenkel, Yi Huang, Athanase Papadopoulos, Dennis Sullivan and Anton Zeitlin for helpful comments and questions.} }}}

\setcounter{footnote}{0}
\section*{Introduction}

In 1939 Wigner   \cite{wigner} described the action of orientation-preserving three-dimensional hyperbolic isometries
as Hermitean conjugation by the complex special linear group.
(See \cite{quinn} and its references for the detailed history going back to Hamilton.)
Specifically, Wigner's action of
$g\in {\rm SL}(2,{\mathbb C})$
on Minkowski space ${\mathbb R}^{3,1}\approx {\mathbb R}^2\times{\mathbb C}\ni(x_1,x_2,x)$ with its quadratic form $x_1x_2-x\bar x$, is
given by 
$A\mapsto g^\dagger A g$, where
$A=\begin{psmallmatrix}x_1&\bar{x}\\x&x_2\\\end{psmallmatrix}$
is a Hermitean matrix and
$ g^{\dagger}=\bar{g}^t$ is the conjugate transpose.
In two dimensions, $g\in {\rm SL}(2,{\mathbb R})$ likewise acts by hyperbolic isometry on $A=\begin{psmallmatrix}
x_1&x\\x&x_2\\\end{psmallmatrix}\in{\mathbb R}^{2,1}$
via $A\mapsto g^t A g$ preserving $x_1x_2-x^2$.

This confluence of two- and three-dimensional algebra and geometry is typically
regarded as a {\sl  low-dimensional accident}, both in geometry and in Lie theory.
We prove here that this accident extends to the super case as well. 
Namely, we 
extend Wigner's action to an action of the complex orthosymplectic group
${\rm OSp}_{\mathbb C}(1{\rm |}2)$ on the super Minkowski space ${\mathbb R}^{3,1|4}$
with its bosonic Wigner coordinates $x_1,x_2$ and $x$ as before plus two complex (or {\it Dirac} as opposed to real {\it Majorana}) fermions $\phi,\psi$ endowed with the quadratic form $x_1x_2-x\bar{x}+\phi\psi +\bar{\phi}\bar{\psi}$.

The discovery of this
super geometry on ${\mathbb R}^{3,1|4}$ which is invariant under a modied
Hermitean conjugation is our main achievement here, a finding of clear mathematical but questionable
physical significance.  Earlier work \cite{others} studied several super metrics on the ${\mathcal  N}=1$ (that is, with
$2{\mathcal N}=2$ fermions, as we have here) super upper half space, but these differ from our metric and are
unrelated to Hermitean conjugation; see the closing remarks in Section \ref{sec:close} for a further discussion.

For $g=\begin{psmallmatrix}a&b&\alpha\\c&d&\beta\\\gamma&\delta&f\\\end{psmallmatrix}\in{\rm OSp}_{\mathbb R}(1{\rm |}2)$
and $A=\begin{psmallmatrix}x_1&x&\phi\\x&x_2&\psi\\
-\phi&-\psi&0\\\end{psmallmatrix}\in {\mathbb R}^{3,1{\rm |}2}$ in the real case 
studied in \cite{cosine,defect,PZ}, the action is given by
$A\mapsto g^{\rm st} A g$, where $g^{\rm st}=\begin{psmallmatrix}a&c&\gamma\\b&d&\delta\\-\alpha&-\beta&f\\\end{psmallmatrix}$
is the super transpose.
A notable point for the sequel is that the bottom-right entry ``0'' in $A$ is preserved under this action in the real case. 
The  invariant form $x_1x_2-x^2+2\phi\psi$ is given by a multiple of the Killing form on the real Lie algebra in this case\footnote{To see this, note that by  definition in the orthosymplectic group, we have $g^{\rm st} Jg=J$, where
$J=\begin{psmallmatrix} 0&-1&0\\1&0&0\\0&0&1\\\end{psmallmatrix}$.  It follows that
$g^{-1}=J^{-1}g^{\rm st} J$, so $g^{-1}Bg=J^{-1}g^{\rm st} J g$.  The adjoint action $B\mapsto g^{-1}Bg$ of course leaves invariant the Killing form and corresponds to our action $A\mapsto g^{\rm st} A g$ for $A=JB$.  However, note that ${\rm OSp}_{\mathbb R}(1{\rm |}2)$ is not invariant under left multiplication by $J$, which simply gives a useful change of coordinates
on the Lie algebra.}.

Berezin has studied the compact real form 
of ${\rm OSp}_{\mathbb C}(1{\rm |}2)$ in \cite{ber}, and a classical construction then
provides a natural Hermitean pairing on the complex Lie algebra derived from
its Killing form.
This Hermitean pairing is invariant under Hermitean conjugation
only by the real form, not by the full complex Lie super group.  The upshot is that classical super Lie theory provides no natural geometry which is invariant under Hermitean conjugation, as far as we can see.

The quadratic form $Q$ introduced here is distinguished by the fact that it is
invariant under the following modified Hermitean conjugation action of ${\rm OSp}_{\mathbb C}(1{\rm |}2)$  on
${\mathbb R}^{3,1{\rm |}2}$, which extends Wigner's action:

\bigskip

\noindent{\bf Main Theorem}.~{\it Given $g=\begin{psmallmatrix}a&b&\alpha\\c&d&\beta\\\gamma&\delta&f\\\end{psmallmatrix}\in{\rm OSp}_{\mathbb C}(1{\rm |}2)$ together with a point $(x_1,x_2,x{\rm |}\phi,\psi)\in{\mathbb R}^{3,1|4}$, define the auxiliary parameter
$$\vartheta=-{{f\bar f}\over 2}[X+(Y-\bar Y)],$$
where $
X=x_1\bar\alpha\alpha+x_2\bar\beta\beta+x\bar\beta\alpha +\bar x\bar\alpha\beta$ and $
Y=\bar\alpha\phi+\bar\beta\psi,
$
and define
$A=\begin{psmallmatrix}x_1&\bar z&\phi\\z&x_2&\psi\\
-\bar\phi& -\bar\psi&\vartheta&\end{psmallmatrix}$.
Then the $g\in{\rm OSp}_{\mathbb C}(1{\rm |}2)$-action 
$A\mapsto g^\dagger A g$ leaves invariant the quadratic form
$$Q(x_1,x_2,z{\rm |}\phi,\psi)=x_1x_2-x\bar{x}+\phi\psi +\bar{\phi}\bar{\psi},$$
where $g^\dagger=\bar g^{\rm st}$ is the conjugate super transpose.}

\medskip

The transformation on super Minkowski space is thus described by Hermitean conjugation,
depending upon an {\it auxiliary parameter} $\vartheta$, followed by projection onto the subspace spanned by $x_1,x_2,x, \phi,\psi$.
We shall write $g.(x_1,x_2,x|\phi,\psi)\in{\mathbb R}^{3,1|4}$ for this action ${\rm OSp}_{\mathbb C}(1|2)\acts {\mathbb R}^{3,1|4}$.  
The proof of the theorem amounts to recognizing symmetries in the formulas
and devolves to computing the roots of a certain quadratic polynomial in $\vartheta$ whose coefficients lie in a Grassmann algebra.  

This seems to be a new
insight.  However, we would be pleased to learn that the auxiliary parameter, a similar extended Hermitean action or invariant quadratic form have arisen in other contexts.

 Interest in the action of ${\rm OSp}_{\mathbb C}(1{\rm |}2)$ on super Minkowski space derives from a desire to extend to three real dimensions aspects of the earlier computations
in two dimensions, where basics on super geodesics and triangles were introduced in \cite{cosine}, and the non-trivial fermionic correction to the Angle Defect Theorem was computed in \cite{defect}.

Section \ref{sec:back} quickly but completely explains the orthosymplectic Lie super groups
${\rm OSp}_{\mathbb R}(1|2)$ and ${\rm OSp}_{\mathbb C}(1|2)$ and their Lie super algebras, as well as a complex structure on a subspace of Minkowski space, which seems not to be well known in this context. The Main Theorem is proved in the Appendix and discussed in Section \ref{sec:wigner}.

In either the real or complex case armed with the Main Theorem, the hyperbolic space $\hyperp$ representing 
$\hyperp^2$ or $\hyperp^3$ is characterized as the locus in super Minkowski space where the invariant quadratic form takes value unity equipped with its induced metric.  In fact, the complex 3d case with Dirac fermions restricts to the real 2d case with Majorana fermions, reflecting
the equivariant embedding $\hyperp^2\subset\hyperp^3$.

Simply this common structure of $\hyperp$ allows for a unified treatment of super geodesics and triangles in either case, including the usual Hyperbolic Law of Cosines for super triangles, as reformulated from \cite{cosine,defect} here in Section \ref{sec:lines} for 3d.  (Presumably this reflects a still more general context for certain of our arguments in super projective geometry.)

The ${\rm OSp}_{\mathbb C}(1|2)$-invariant three-form on
$\hyperp^3$ whose body is the hyperbolic volume form is computed in Section \ref{sec:vol}, and a 
primitive two-form is derived for it.  Super triangles and tetrahedra are discussed in Section \ref{tris} including divergence of the volume of an ideal tetrahedron.
A final Section \ref{sec:close} of closing remarks in particular explains the failure of the Schl\"afli formula for super tetrahedra as well as ideas for computing its non-vanishing fermionic correction.

Though this note is self-contained, we refer the reader to: \cite{cosine,defect} for background on
the real orthosymplectic group and its action on the super Minkowski space
${\mathbb R}^{2,1|2}$; for applications beyond the scope and rather independent from much of what is discussed here 
for ${\mathcal N}=1$ super Teichm\"uller space to \cite{PUN,PZ}; and for ${\mathcal N}=2$ to
\cite {IPZ}.

\section{Orthosymplectic Groups and Super Minkowski Spaces}\label{sec:back}

In this section, we include basic information concerning the real and complex orthosymplectic groups ${\rm OSp_{\mathbb R}(1|2)}$ and ${\rm OSp_{\mathbb C}(1|2)}$, which are among the simplest of Lie super groups and whose respective bodies are the classical special linear groups ${\rm SL}(2,{\mathbb R})$ and ${\rm SL}(2,{\mathbb C})$.  We refer the interested reader to \cite{Kac,Manin}
for more information about general Lie super algebras and super groups
and to \cite{PZ} for details about OSp$_{\mathbb R}$(1$\vert$2).

\subsection{Super Numbers}

Let $\hat{{\mathbb R}}=\hat{\mathbb R}[0]\otimes\hat{\mathbb R}[1]$ 
and $\hat{{\mathbb C}}=\hat{\mathbb C}[0]\otimes\hat{\mathbb C}[1]$ be the respective ${\mathbb Z}$/2-graded
rings of polynomials over ${\mathbb R}$ and over ${\mathbb C}$ with one central generator
$1\in {\mathbb R}\subseteq \hat{\mathbb R}[0]\subseteq \hat{\mathbb C}[0]$ of degree zero and infinitely many anti-commuting generators $\theta_1,\theta_2,\ldots$ of degree one.
An arbitrary $a\in\hat{\mathbb R}$ or $\hat{\mathbb C}$ can be written uniquely as
$$a=a_{\#}+\sum_i a_i\theta_i+\sum_{i<j} a_{ij} \theta_i \theta_j++\sum_{i<j<k} a_{ijk} \theta_i \theta_j \theta_k\cdots,$$
where the coefficients $a_\#, a_i,a_{ij},a_{ijk}$ are respectively real or complex.  
The  term $a_\#$ of degree zero is called the {\it body} of the {\it super number} $a\in \hat{\mathbb R}$.  
A super number is invertible if and only if its body is non-zero.

If $a\in\hat{\mathbb C}[0]$, then it is said to be an {\it even} super number or  {\it boson}, while if $a\in\hat{\mathbb C}[1]$, then it is said to  be an {\it odd} super number or {\it fermion}. We adopt the notation throughout that fermions are denoted by ordinary lower-case Greek letters, and typically conversely.

The usual order relation $\leq$ on ${\mathbb R}$ induces one on ${\hat{\mathbb R}}$ with $a\leq b$, for $a,b\in\hat{\mathbb R}$, if and only if $a_\#\leq b_\#$.  Likewise if $a_\#\neq 0$, then the sign $sign(a)$ is defined to be the sign of $a_\#$ and the absolute value is $|a|=sign(a)a$.  Other obvious extensions are implicit.

\subsection{The Lie super groups}

An element $g$ in the {\it orthosymplectic group} can be represented by the 3-by-3 matrix $g=\begin{psmallmatrix}
a&b&\alpha\\c&d&\beta\\
\gamma&\delta&f\\\end{psmallmatrix}$, where $a,b,c,d,f$ are even and $\alpha,\beta,\gamma,\delta$ are odd and all entries are real or complex in the respective cases ${\rm OSp}_{\mathbb R}(1{\rm |}2)$
and ${\rm OSp}_{\mathbb C}(1{\rm |}2)$, 
with multiplication in either case defined by
$$\biggl ( \begin{smallmatrix}
a_1&b_1&\alpha_1\\c_1&d_1&\beta_1\\\gamma_1&\delta_1&f_1\\
\end{smallmatrix}\biggr )
\biggl ( \begin{smallmatrix}
a_2&b_2&\alpha_2\\c_2&d_2&\beta_2\\\gamma_2&\delta_2&f_2\\
\end{smallmatrix}\biggr )=
\biggl (\begin{smallmatrix}
a_1a_2+b_1c_2-\alpha_1\gamma_2&a_1b_2+b_1d_2-\alpha_1\delta_2&a_1\alpha_2+b_1\beta_2+\alpha_1 f_2\\
c_1a_2+d_1c_2-\beta_1\gamma_2&c_1b_2+d_1d_2-\beta_1\delta_2&c_1\alpha_2+d_1\beta_2+\beta_1 f_2\\
\gamma_1 a_2+\delta_1 c_2+\delta_1\gamma_2&\gamma_1 b_2+\delta_1 d_2+f_1 \delta_2&-\gamma_1\alpha_2-\delta_1\beta_2+f_1f_2\\
\end{smallmatrix}\biggr ),$$
and satisfying the two conditions

$\bullet$~the {\it super determinant} or {\it Berezinian} of $g$ is unity, namely,
$$
{\rm sdet}~g~=~f^{-1}
~\det\left[
\begin{pmatrix}
a&b\\ c&d
\end{pmatrix}
+f^{-1}
\begin{pmatrix}
\alpha\gamma&\alpha\delta\\ \beta\gamma&\beta\delta
\end{pmatrix}
\right]=1;
$$ 

$\bullet$~$g$ is {\it orthosymplectic}, namely
$$g^{st}Jg=J,$$
where
$
J= \begin{psmallmatrix}
~0 & -1 &0 \\
\hskip .8ex1 &\hskip1.2ex 0 &0 \\
~0 &\hskip1.2ex 0 & 1 \end{psmallmatrix}
$
and the {\it super transpose} $g^{st}=\begin{psmallmatrix}
~a & ~c & \gamma \\
~b & ~d & \delta \\
-\alpha & -\beta & f \\
\end{psmallmatrix}.$

Routine computations show that these two conditions are
equivalent to the defining relations
\begin{eqnarray}\nonumber
&&\alpha=b\gamma-a\delta, \quad \beta =d\gamma-c\delta, \quad \hskip 2.8ex f=1+\alpha\beta,\nonumber\\
&&\gamma=a\beta-c\alpha, \quad \delta=b\beta-d\alpha, \quad f^{-1}=ad-bc.\nonumber
\end{eqnarray}
for the orthosymplectic groups, so in particular $\alpha\beta=\gamma\delta$.

As in the real case, the subgroup 
$$\widehat{\rm SL}(2,{\mathbb C})=
\biggl \{\begin{pmatrix}a&b&0\\c&d&0\\0&0&1\\\end{pmatrix}\in{\rm OSp}_{\mathbb C}(1{\rm |}2):
a,b,c,d\in\hat{\mathbb C}[0]\biggr\}$$
is intermediary between the classical ${\rm SL}(2,{\mathbb C})$
and the full ${\rm OSp}_{\mathbb C}(1{\rm |}2)$.  Another useful subspace, which is not a subgroup of
${\rm OSp}_{\mathbb C}(1|2)$, is defined by setting $$u(\alpha,\beta)=\begin{pmatrix}1-{{\alpha\beta}\over 2}&0&\alpha\\0&1-{{\alpha\beta}\over 2}&\beta\\\beta&-\alpha&1+\alpha\beta\\\end{pmatrix}\in{\rm OSp}_{\mathbb C}(1|2),
~{\rm for}~\alpha,\beta\in\hat{\mathbb C}[1].$$

\begin{lemma}\label{factor}
Any element $g\in{\rm OSp}_{\mathbb C}(1|2)$ can be written uniquely as a product
$$g=\begin{psmallmatrix}a&b&0\\c&d&0\\0&0&1\\\end{psmallmatrix}\,u(\alpha,\beta)=u(a\alpha+b\beta,c\alpha+d\beta)\, \begin{psmallmatrix}a&b&0\\c&d&0\\0&0&1\\\end{psmallmatrix},$$
for some $\begin{psmallmatrix}a&b&0\\c&d&0\\0&0&1\\\end{psmallmatrix}\in \widehat{\rm SL}(2,{\mathbb C})$ and fermions $\alpha,\beta\in\hat{\mathbb C}[1]$.\hfill$\mathlarger{\mathlarger{\mathlarger{\mathlarger {\square}}}}$
\end{lemma}

\subsection{The Lie super algebras}
Respective Lie super algebras ${\rm osp}_{\mathbb R}(1|2)$ and ${\rm osp}_{\mathbb C}(1|2)$ of
${\rm OSp}_{\mathbb R}(1{\rm |}2)$
and ${\rm OSp}_{\mathbb C}(1{\rm |}2)$ are characterized by the functional equation 
$A+JA^{ st}=0$
and have three even $h, X_{\pm}$ and two odd generators $v_{\pm}$ satisfying the super commutation relations
$$
[h,v_{\pm}]=\pm v_{\pm}, \quad[v_{\pm},v_{\pm}]=\mp 2X_{\pm},\quad  
[v_+,v_{-}]=h,
$$
where the bracket denotes the relevant commutator or anti-commutator.
These are represented by $(2|1)\times (2|1)$ supermatrices given by
\begin{eqnarray}
v_{+}=\left( \begin{array}{ccc}
0 & \hskip1ex0 & 1 \\
0 & \hskip1ex0 & 0 \\
0 & -1 & 0 \end{array} \right), \quad 
v_{-}=\left( \begin{array}{ccc}
0 & 0 & 0 \\
0 & 0 & 1 \\
1 & 0 & 0 \end{array} \right), \quad 
h=\left( \begin{array}{ccc}
1 & \hskip1ex0 & 0 \\
0 & -1 & 0 \\
0 & \hskip1ex0 & 0 \end{array} \right) .\nonumber
\end{eqnarray}

The seemingly exotic signs in products of matrices in the super group in the previous section
assure the usual products and signs in the graded super algebra in this section.
Furthermore, even though left multiplcation by $J$ does not preserve the orthosymplectic super Lie algebras, we shall find it anyway useful to consider the change of coordinates $A=JB$ for $B\in {\rm osp}_{\mathbb R}(1{\rm |}2)$ or $B\in{\rm osp}_{\mathbb C}(1{\rm |}2)$.

\subsection{Super Minkowski Space}
As in \cite{cosine,defect}, super real Minkowski space of dimension (2,1$|$2) is
$${\mathbb R}^{2,1|2}=\{ (x_1,x_2,x~\vert~\phi,\psi)\in\hat{\mathbb R}^5:x_1,x_2,x\in\hat{\mathbb R}[0]~{\rm and}~ \phi,\psi\in\hat{\mathbb R}[1]\}$$
endowed wih the quadratic form $x_1x_2-x^2+2\phi\psi$, which is essentially the Killing form of ${\rm OSp}_{\mathbb R}(1|2)$.  In this paper, we shall consider the analogous super real Minkowski space 
$${\mathbb R}^{3,1|4}=\{ X=(x_1,x_2,x~\vert~\phi,\psi)\in\hat{\mathbb C}^5:x_1,x_2\in\hat{\mathbb R}[0], x\in\hat{\mathbb C}[0], \phi,\psi\in\hat{\mathbb C}[1]\}$$
of dimension (3,1$|$4), where one bosonic and both fermionic coordinates have been 
complexified from the real case, as will be further explained.

Writing a vector as a matrix for real $A=JB=\begin{psmallmatrix}
x_1&x&\phi\\x&x_2&\psi\\-\phi&-\psi&0\\
\end{psmallmatrix}$, we have
$$B=\psi v_+-\phi v_-+xh+x_2X_+-x_1X_-,$$
and setting $ \hat A=\begin{psmallmatrix}
x_1&\bar z&\phi\\z&x_2&\psi\\-\bar\phi&-\bar\psi&0\\
\end{psmallmatrix}$ with three coordinates complexified, we
define the further 
\begin{eqnarray}
\hat v_{+}=\left( \begin{array}{ccc}
0 & \hskip1ex0 & 1 \\
0 & \hskip1ex0 & 0 \\
0 & 1 & 0 \end{array} \right), \quad 
\hat v_{-}=\left( \begin{array}{ccc}
0 & 0 & 0 \\
0 & 0 & 1 \\
-1 & 0 & 0 \end{array} \right), \quad 
\hat h=\left( \begin{array}{ccc}
1 & \hskip1ex0 & 0 \\
0 & 1 & 0 \\
0 & \hskip1ex0 & 0 \end{array} \right) ,\nonumber
\end{eqnarray}
so that 
$$\begin{aligned}
\hat B&=J^{-1}\hat A\\
&=({\rm Re}\,\psi~ v_++i{\rm ~Im}\,\psi~ \hat v_+)
-({\rm Re}\,\phi~ v_-+i{\rm ~Im}\,\phi~ \hat v_-)\\
&+({\rm Re}\,~z~ h~+i{\rm ~Im}\, z~ \hat h)+x_2X_+-x_1X_-.\\
\end{aligned}$$ 
It follows that
$$\begin{aligned}
2\psi{\partial\over{\partial\alpha}}&=\psi(v_++i\hat v_+)+\bar\psi(v_+-i\hat v_+),\\
-2\phi{\partial\over{\partial\beta}}&=\phi(v_-+i\hat v_-)+\bar\phi(v_--i\hat v_-),\\
\end{aligned}$$
explaining the implicit complex structure on this sub-bundle 
in ${\mathbb R}^{3,1|4}$ of the real bundle
underlying the complex Lie super algebra.




\section {Super Wigner Representation}\label{sec:wigner}

In the (2,1$|$2)-dimensional case with $g=\begin{psmallmatrix}a&b&\alpha\\c&d&\beta\\\gamma&\delta&f\\\end{psmallmatrix}\in{\rm OSp}_{\mathbb R}(1{\rm |}2)$
and regarding $(x_1,x_2,y|\phi,\psi)\in {\mathbb R}^{3,1{\rm |}2}$ as a matrix
$A=\begin{psmallmatrix}x_1&y&\phi\\y&x_2&\psi\\
-\phi&-\psi&0\\\end{psmallmatrix}\in {\mathbb R}^{3,1{\rm |}2}$, the real orthosymplectic action is given by
$A\mapsto g^{\rm st} A g$, where the super transpose is $g^{\rm st}=\begin{psmallmatrix}a&c&\gamma\\b&d&\delta\\-\alpha&-\beta&f\\\end{psmallmatrix}$
.

Again in the complex case given
$g=\begin{psmallmatrix}a&b&\alpha\\c&d&\beta\\\gamma&\delta&f\\\end{psmallmatrix}\in{\rm OSp}_{\mathbb C}(1{\rm |}2)$,
we write a matrix
$A=\begin{psmallmatrix}x_1&\bar x&\phi\\x&x_2&\psi\\
-\bar\phi& -\bar\psi&\vartheta&\end{psmallmatrix}$
for  a point $(x_1,x_2,x|\phi,\psi)\in{\mathbb R}^{3,1|4}$ now with 
$$\vartheta=-{{f\bar f}\over 2}[X+(Y-\bar Y)],$$
where $
X=x_1\bar\alpha\alpha+x_2\bar\beta\beta+x\bar\beta\alpha +\bar x\bar\alpha\beta$ and $
Y=\bar\alpha\phi+\bar\beta\psi.
$
The Main Theorem from the Introduction asserts that
the quadratic form $Q(x_1,x_2,x{\rm |}\phi,\psi)=x_1x_2-x\bar{x}+\phi\psi +\bar{\phi}\bar{\psi}$
is invariant under the action $A\mapsto g^\dagger A g=\begin{psmallmatrix}x_1'&\bar x'&\phi'\\x'&x'_2&\psi'\\
-\bar\phi'& -\bar\psi'&\vartheta'&\end{psmallmatrix}$.  We may sometimes suppress $\vartheta$ entirely and write simply
$g.(x_1,x_2,x|\phi,\psi)=(x_1',x_2',x'|\phi',\psi')$.

The proof is a direct but lengthy computation at certain points concise and rewarding and at others
routine and tedious, which is relegated to the Appendix.   As a corollary to this proof, we have

\begin{corollary}\label{corgood}
Suppose that $g.(x_1,x_2,x|\phi,\psi)=(x_1',x_2',x|'\phi',\psi')$ with auxiliary parameter $\vartheta$
in the notation above.  Then
we have the identities
$$\begin{aligned}
-2\vartheta&=(1+\bar\alpha\bar\beta)(\bar\phi\alpha+\bar\psi\beta)-
(1+\alpha\beta)(\phi\bar\alpha+\psi\bar\beta)\\&+\begin{pmatrix}\bar\alpha &\bar\beta\\\end{pmatrix}\begin{pmatrix}
x_1&\bar x\\x&x_2\\\end{pmatrix}\begin{pmatrix}\alpha\\\beta\\\end{pmatrix}
,\\
\begin{pmatrix}\bar a&\bar c\\\bar b&\bar d\\
\end{pmatrix}^{-1}
\begin{pmatrix}\phi'\\\psi'\\\end{pmatrix}=&
\biggl(1+\alpha\beta+{{\bar\alpha\bar\beta}\over 2}(1+2\alpha\beta)\biggr)
\begin{pmatrix}\phi\\\psi\\\end{pmatrix}\\
&+\begin{pmatrix}
x_1(1+{{\bar\alpha\bar\beta}\over 2})&\bar x(1+{{\bar\alpha\bar\beta}\over 2})\\
x(1-{{\bar\alpha\bar\beta}\over 2})&x_2(1-{{\bar\alpha\bar\beta}\over 2})
\end{pmatrix}\begin{pmatrix}\alpha\\\beta\\\end{pmatrix},\\
\begin{pmatrix}\bar a&\bar c\\\bar b&\bar d\\
\end{pmatrix}^{-1}\begin{pmatrix}x_1'&\bar x'\\x'&x_2'\end{pmatrix}&
\begin{pmatrix} a&b\\c&d\\
\end{pmatrix}^{-1}=\biggl(1-{{\alpha\beta\bar\alpha\bar\beta}\over 2}\biggr)\begin{pmatrix}x_1&\bar x\\x&x_2\\\end{pmatrix}\\
+\biggl (1+{{\alpha\beta}\over 2}\biggr)&\begin{pmatrix}\bar\beta\\-\bar\alpha\\\end{pmatrix}\begin{pmatrix}\bar\phi&\bar\psi\end{pmatrix}
-\biggl (1+{{\bar\alpha\bar\beta}\over 2}\biggr)\begin{pmatrix}\phi\\\psi\\\end{pmatrix}\begin{pmatrix}\beta&-\alpha\\\end{pmatrix}.
\end{aligned}$$
\end{corollary}

\begin{proof} The formula for $\vartheta$ follows immediately from the definitions 
using the usual identity
 $\begin{psmallmatrix} \gamma\\\delta\\\end{psmallmatrix}=\begin{psmallmatrix} a&c\\ b&d\\\end{psmallmatrix}\begin{psmallmatrix} \beta\\-\alpha\\\end{psmallmatrix}$
 for $g\in{\rm OSp}(1|2)$. It likewise follows immediately from the definitions that
$$\begin{aligned}
\begin{pmatrix}\bar a&\bar c\\\bar b&\bar d\\
\end{pmatrix}^{-1}\begin{pmatrix}\phi'\\\psi'\end{pmatrix}&= (1+\alpha\beta)\begin{pmatrix}\phi\\\psi\\\end{pmatrix}
~+~\begin{pmatrix}x_1&\bar x\\x&x_2\end{pmatrix}\begin{pmatrix}\alpha\\\beta\\\end{pmatrix}\\
&+\biggl\{
\bar\phi \alpha+\bar\psi\beta+(1+\alpha\beta)\vartheta\biggr\}\begin{pmatrix}\bar\beta\\-\bar\alpha\\\end{pmatrix},\\
\begin{pmatrix}\bar a&\bar c\\\bar b&\bar d\\
\end{pmatrix}^{-1}\begin{pmatrix}x_1'&\bar x'\\x'&x_2'\end{pmatrix}&
\begin{pmatrix} a&b\\c&d\\
\end{pmatrix}^{-1}=
\begin{pmatrix}x_1&\bar x\\x&x_2\end{pmatrix}-\vartheta\begin{pmatrix}\bar\beta\\-\bar\alpha\end{pmatrix}
\begin{pmatrix}\beta&-\alpha\end{pmatrix}\\
&\hskip 10ex+\begin{pmatrix}\bar\beta\\-\bar\alpha\end{pmatrix}\begin{pmatrix}\bar\phi&\bar\psi\end{pmatrix}
-\begin{pmatrix}\phi\\\psi\end{pmatrix}\begin{pmatrix}\beta&-\alpha\end{pmatrix}.\\
\end{aligned}$$
Substitution of the expression for $\vartheta$ into these formulae provides the asserted identities by routine computation.
\end{proof}

Notice that this action is linear in the super algebra variables and has higher-order
fermionic corrections to being quadratic in the super group variables.  Furthermore,
$\vartheta$ vanishes for $g\in\widehat{\rm SL}(2,{\mathbb C})$.  Moreover, the bottom-right entry
$\vartheta '$ of
of $g^\dagger \hat A g$ is not invariant, but rather $\vartheta '=-2\vartheta$, for which we have neither explanation nor application.

Let 
$$\langle \cdot\,, \cdot\,\rangle: {\mathbb R}^{3,1|4}\times
{\mathbb R}^{3,1|4}\to\hat{\mathbb R}[0]$$
denote the bilinear form associated with the quadratic form $Q$ in the Main Theorem, i.e.,
$$\begin{aligned}
&2\langle(x_1,x_2,x|\phi,\psi),(y_1,y_2,y|\xi,\eta)
\rangle\\
&=
x_1y_2+x_2y_1-(x\bar y+\bar x y)+\phi\eta+\xi\psi+\bar\phi\bar\eta+\bar\xi\bar\psi.
\end{aligned}
$$Note that
this restricts for real $x$, Majorana fermions, and the real orthosymplectic group action
to the invariant real form, and so we use the same notation for the inner product
in either case.

Define our models
$$\begin{aligned}
\hyperp^2&=\{ X=(x_1,x_2,x|\phi,\psi)\in{\mathbb R}^{2,1|2}:\langle X,X\rangle=1~{\rm and}~x_1+x_2>0\},\\
\hyperp^3&=\{ X=(x_1,x_2,x|\phi,\psi)\in{\mathbb R}^{3,1|4}:\langle X,X\rangle=1~{\rm and}~x_1+x_2>0\}\\
\end{aligned}$$
for the super hyperbolic spaces, with the metric induced by the bilinear form,
let
$$\begin{aligned}
\hyperl^2&=\{ X=(x_1,x_2,x|\phi,\psi)\in{\mathbb R}^{2,1|2}:\langle X,X\rangle=-1\},\\
\hyperl^3&=\{ X=(x_1,x_2,x|\phi,\psi)\in{\mathbb R}^{3,1|4}:\langle X,X\rangle=-1\}\\
\end{aligned}$$
denote the super hyperboloids of one sheet, and let
$$\begin{aligned}
{\mathcal L} _+^2&=\{ X=(x_1,x_2,x|\phi,\psi)\in{\mathbb R}^{2,1|2}:\langle X,X\rangle=0~{\rm and}~x_1+x_2>0\},\\
{\mathcal L} _+^3&=\{ X=(x_1,x_2,x|\phi,\psi)\in{\mathbb R}^{3,1|4}:\langle X,X\rangle=0~{\rm and}~x_1+x_2>0\}\\
\end{aligned}$$
denote the super positive light-cones in each case.  

When the dimension 2 or 3 is fixed or immaterial, we may write simply $\hyperp,\hyperl$, or ${\mathcal L}_+$
or even ${\rm OSp}(1|2)$.
It is also useful to have notation for the underlying bodies, to wit
$\widecheck\hyperp\subseteq \hyperp, \widecheck\hyperl\subseteq\hyperl$, and $\widecheck{\mathcal L}_+\subseteq{\mathcal L}_+$.

\section{Super geodesics}\label{sec:lines}

This section is devoted to the proof of the following theorem, which is
adapted from \cite{defect} and applies equally well in the earlier real and the current complex
cases, as it depends only formally on the existence of the inner product and hyperboloid.

\begin{theorem} \label{thm:lines}
Given distinct $P,Q\in\hyperp$, there is a unique super geodesic containing them.
In fact, define
$d=\langle P,Q\rangle$ and set 
$\ell=\sqrt{{d+1}\over{d-1}}$.  Then the distance $D$ between $P$ and $Q$ is given by
${\rm cosh}\, D=d$ with ${\rm exp}\, D={{\ell+1}\over{l-1}}$.
The super segment $\overline {PQ}$ from $P$ to $Q$ is parametrized by
$$X(t)={{tE+t^{-1}F}\over 2},~{\rm for}~1\leq t\leq {{\ell+1}\over{\ell-1}},$$
where $t=~{\rm exp}\, s$, with $s$ equal arc length, and 
$$\begin{aligned}E&=M\{(1-\ell)P+(1+\ell)Q\},\\
F&=N\{(1+\ell)P+(1-\ell)Q\}\end{aligned}$$
lie in ${\mathcal L}_+$, and we take $M={{\ell-1}\over {2\ell}}, N={{\ell+1}\over{2\ell}}$
so that $\langle E,F\rangle=2$.
\end{theorem}

 \begin{proof} We first claim that
the general form of a super geodesic  
${ X}(s)$ in $\hyperp$ parametrized by arc length $s$ is given by
$$
{X}=\cosh s\,{ U}+\sinh s\,{ V},
$$
for some ${ U}\in\hyperp$, ${ V}\in\hyperl$ with $\langle{ U}, { V}\rangle =0$. This follows directly from the variational principle applied to the functional 
$$
\int \Big({ |\langle \dot{ X}, \dot{ X}\rangle|}+\lambda( \langle { X}, { X}\rangle-1)\Big)dt ,
$$
where dot stands for ${d\over{ds}}$ along the curve,
with Euler-Lagrange equations 
$
{\ddot{ X}}=\lambda { X},~~ \langle { X}, { X}\rangle=1 \label{geod}
$,
and with $s$ chosen so that $\langle\dot{ X}, \dot{ X}\rangle=1$.

For any ${ U}\in\hyperp$ and ${ V}\in\hyperl$ with $\langle { U},{ V}\rangle =0$, let 
$$L_{{ U},{ V}}=\{{\rm cosh}\, s~{ U}+{\rm sinh}\, s~{ V}:s\in{\mathbb R}\}$$ denote the corresponding super geodesic, whose asymptotes are given by the rays in ${\mathcal L}_+$
containing the vectors
$ E=U+V$, $F=U-V$.
Conversely,  points ${ E},{ F}\in{\mathcal L}_+$ with $\langle { E},{ F}\rangle =2$ determine a unique corresponding super geodesic, where we have ${ U}={1\over 2}({ E}+{ F})\in \hyperp$ and ${ V}={1\over 2}({ E}-{ F})\in\hyperl$.

To prove existence of a super geodesic containing $P$ and $Q$, we have $\langle{ P},{ Q}\rangle~>1$
for distinct $P,Q\in\hyperp$,  so we may define
$$\ell=\sqrt{{\langle{ P},{ Q}\rangle+1}\over{\langle{ P},{ Q}\rangle-1}},~~
{\rm whence}~\langle{ P},{ Q}\rangle={{\ell^2+1}\over{\ell^2-1}}.$$  
The identity
${\rm cosh}^{-1}t={\rm log}_e(t+\sqrt{t^2-1})$ therefore gives 
$${\rm exp}\,D={\rm exp}\,{{\rm cosh}^{-1}\langle{ P},{ Q}\rangle}=\langle{ P},{ Q}\rangle+\sqrt{\langle{ P},{ Q}\rangle^2-1}={{\ell+1}\over{\ell -1}}.$$

We exhibit $L_{{ U},{ V}}=L_{{{{ E}+{ F}}\over 2},{{{ E}-{ F}}\over 2}}$ containing ${ P},{ Q}$ parametrized as
$${ X}(s)={\rm cosh}\, s~{ U}+{\rm sinh}\, s~{ V}={1\over 2}[{\rm exp}(s)~{ E}+{\rm exp}(-s)~{ F}],$$
where
$$\begin{aligned}
{ E}&={{\ell-1}\over{2\ell}}[(1-\ell){ P}+(1+\ell){ Q}],\\
{ F}&={{\ell+1}\over{2\ell}}[(1+\ell){ P}+(1-\ell){ Q}].\\
\end{aligned}$$
Direct computation confirms that ${ E},{ F}\in{\mathcal L}_+$ with $\langle{ E},{ F}\rangle=2$ and
$${ X}(0)={ P}~{\rm and}~{ X}({\rm cosh}^{-1}\langle{ P},{ Q}\rangle)={ Q},$$
thus establishing existence of as well as that ${\rm cosh}\,D=\langle { P},{ Q}\rangle$. 
Uniqueness amounts to invertibility of 
$\begin{psmallmatrix}{\rm cosh}\, s\,I&{\rm sinh}\, s\,I\\{\rm cosh}\, t\,I&{\rm sinh}\, t\,I\end{psmallmatrix},$
where $I$ is the 5-by-5 identity matrix.\end{proof}

It follows that any three non-collinear points of $\hyperp$ define a {\it super triangle}, namely, three super geodesic segments with disjoint interiors meeting pairwise
at the given points.

\begin{corollary}\label{cor:tgt}
Given distinct ${ P},{ Q}\in\hyperp$, the unit tangent vector at ${ P}$ to the super geodesic segment from ${ P}$ to ${ Q}$ is given by
 ${{ Q}-{ P}\,\langle { P},{ Q}\rangle}\over{\sqrt{\langle { P},{ Q}\rangle^2-1}}$.~\hfill \qedsymbol
 \end{corollary}

\begin{corollary} The usual Hyperbolic Law of Cosines holds for super triangles in 
$\hyperp$.
\end{corollary}
\begin{proof}
This follows directly from Corollary \ref{cor:tgt}
upon computing the cosine of an interior angle via the usual formula
as in \cite{cosine}.
\end{proof}

\begin{corollary}\label{useful}
Let $L=L_{{U},{V}}$ be a super geodesic in $\hyperp$  and define the isotropic ${E}={U}+{V}, {F}={U}-{V}\in{\mathcal L}_+$.
Then we have
$$\begin{aligned}
L&=\{ P\in\hyperp:\langle {P},{E}\rangle~\langle {P},{F}\rangle~=~1\}\, ,\\
&=\biggl\{ {1\over{2\sqrt{xy}}}(x{E}+y{F}):x,y> 0\biggr\}\,.\\
\end{aligned}$$
\end{corollary}

\begin{proof} For the inclusion of $L$ in the first equality, write $${P}={\rm cosh}\, p~{U}~+~{\rm sinh}\, p~{V}\,,$$ so that $\langle {P},{U}\pm{V}\rangle={\rm cosh}\, p\mp{\rm sinh}\,p,$ whence
$$\langle {P},{E}\rangle~\langle {P},{F}\rangle={\rm cosh}^2 p-{\rm sinh}^2 p=1.$$

\bigskip

For the reverse inclusion, suppose that
$$\begin{aligned}
1&=\langle {P},{E}\rangle~\langle {P},{F}\rangle\\
&=\langle {P},{U}\rangle^2-\langle {P},{V}\rangle^2\, \\
\end{aligned}$$
and define ${Q}=\langle {P},{U}\rangle\,{U}-\langle {P},{V}\rangle\,{V}$,
whence
$$\begin{aligned}
\langle{Q},{Q}\rangle&=\biggl\langle \langle {P},{U}\rangle\,{U}-\langle {P},{V}\rangle\,{V}~,~\langle {P},{U}\rangle\,{U}-\langle {P},{V}\rangle\,{V}\biggr\rangle\\
&=\langle{P},{U}\rangle^2-\langle{P},{V}\rangle^2\\
&=1,
\end{aligned}$$
so ${Q}\in\hyperp$.  Moreover by Theorem \ref{thm:lines}, the cosh of the distance between
${P}$ and ${Q}$ is given by
$$\begin{aligned}
\langle{P},{Q}\rangle&=\biggr\langle{P}, \langle {P},{U}\rangle\,{U}-\langle {P},{V}\rangle\,{V}\biggr\rangle\\
&=\langle{P},{U}\rangle^2-\langle{P},{V}\rangle^2,\\
\end{aligned}$$
so the distance between ${P}$ and ${Q}$ is zero, whence ${P}={Q}$,
proving the first identity.

\bigskip

For the second equality, suppose
${Q}=x{E}+y{F}$, for $x,y>0,$ so
$$\langle {Q},{E}\rangle~\langle{Q},{F}\rangle~=~4xy=\langle{Q},{Q}\rangle$$
since $\langle{E},{F}\rangle=2$.  The first equality shows that 
${1\over{2\sqrt{xy}}}{Q}\in L_{{U},{V}}$ and the second that ${1\over{2\sqrt{xy}}}{Q}\in\hyperp$ as required.
\end{proof}

\noindent The significant computational import of this result is that a super geodesic is the projectivization to lie in $\hyperp$ of the convex linear span of vectors lying in its
asymptotes in ${\mathcal L}_+$, just as in the body.

\section{Volume Form}\label{sec:vol}

In this section, we compute a three-form Vol on $\hyperp^3$ which is invariant under
${\rm OSp}_{\mathbb C}(1|2)$ and whose body is the hyperbolic volume form; we shall also find
a natural primitive for it.  This is in analogy to \cite{defect} where
the two-form
$\Omega(x_1,x_2,y|\phi,\psi)={{{\rm d}[(1+\phi\psi)x_2)]\wedge{\rm d}[(1+\phi\psi)x_1]}\over {2[(1+\phi\psi)y)]}}$
is invariant under ${\rm OSp}_{\mathbb R}(1|2)$, and
$\Omega={\rm d}\omega$ for
$\omega(x_1,x_2,y|\phi,\psi)={{y(1+\phi\psi)}\over 2}~{\rm d}{\rm log}_e ( {{x_1}\over{x_2}})$.

Likewise analogously to \cite{cosine,defect} in 2d, the following result allows us to pull-back
3d hyperbolic geometry to super hyperbolic space.

\begin{lemma} Given $(x_1,x_2,x|\phi,\psi)\in{\mathbb H}^3$, 
define
$$\begin{pmatrix}\alpha\\\beta\\\end{pmatrix}=(1+{1\over 2}\bar\phi\bar\psi)
\begin{pmatrix}-x_2&\bar x\\x&-x_1\\\end{pmatrix}\begin{pmatrix}\phi\\\psi\\\end{pmatrix}.$$
Then for $u=u(\alpha,\beta)$ in the notation of Lemma \ref{factor}, we have
$$u.(x_1,x_2,x|\phi,\psi)=[1+{1\over 2}(\phi\psi+\bar\phi\bar\psi)+{3\over 4}\phi\psi\bar\psi\bar\psi]
(x_1,x_2,x|0,0).$$
\end{lemma}

\begin{proof}
 It follows from the definitions that
$$\begin{aligned}
\phi'&=(x_1\alpha+\bar x\beta)-{{\bar\beta}\over 2}(\alpha\bar\phi+\beta\bar\psi)
+(1+\alpha\beta+{1\over 2}\,\alpha\beta\bar\alpha\bar\beta)\phi,\\
\psi'&=(x_2\beta+x\alpha)+{{\bar\alpha}\over 2}(\alpha\bar\phi+\beta\bar\psi)
+(1+\alpha\beta+{1\over 2}\,\alpha\beta\bar\alpha\bar\beta)\psi,\\
\end{aligned}$$
so $\phi'=\psi'=0$ is equivalent to
$$\begin{pmatrix}
x_1-{{\bar\phi\bar\beta}\over 2}&\bar x-{{\bar\psi\bar\beta}\over 2}\\
x+{{\bar\phi\bar\alpha}\over 2}&x_2+{{\bar\psi\bar\alpha}\over 2}\\
\end{pmatrix}\begin{pmatrix}\alpha\\\beta\\\end{pmatrix}=-
(1+\alpha\beta+{1\over 2}\,\alpha\beta\bar\alpha\bar\beta)
\begin{pmatrix}\phi\\\psi\\\end{pmatrix}.$$

Thus, we find that
$$\begin{aligned}
-\begin{pmatrix}\phi\\\psi\\\end{pmatrix}
&=\begin{pmatrix}
x_1-{{\bar\phi\bar\beta}\over 2}&\bar x-{{\bar\psi\bar\beta}\over 2}\\
x+{{\bar\phi\bar\alpha}\over 2}&x_2+{{\bar\psi\bar\alpha}\over 2}\\
\end{pmatrix}
\begin{pmatrix}\alpha\\\beta\\\end{pmatrix}\\
&=\biggl\{\begin{pmatrix}x_1&\bar x\\x&x_2\\\end{pmatrix}
+{1\over 2}\begin{pmatrix}\bar\beta\\-\bar\alpha\end{pmatrix}\begin{pmatrix}\bar\phi&\bar\psi\\\end{pmatrix}\biggr\}
\begin{pmatrix}\alpha\\\beta\\\end{pmatrix},
\end{aligned}$$
and using the constraint $x_1x_2-x\bar x+\phi\psi+\bar\phi\bar\psi$=1, we conclude
$$\begin{aligned}
-\begin{pmatrix}\alpha\\\beta\end{pmatrix}&=
\begin{pmatrix}x_2&-\bar x\\-x&x_1\\\end{pmatrix}\\
&\biggl[
(1+\bar\phi\bar\psi)\begin{pmatrix} \phi\\\psi\end{pmatrix}
+{1\over 2}(1+\phi\psi)(\bar\phi\alpha+\bar\psi\beta)\begin{pmatrix}\bar\beta\\-\bar\alpha\end{pmatrix}
\biggr].\end{aligned}$$
In particular, a short computation again using the constraint 
gives
$$\begin{aligned}
\alpha\beta&=\phi\psi(1+\bar\phi\bar\psi)-{1\over 2}(1+\phi\psi)\bar\phi\bar\psi\alpha\beta\bar\alpha\bar\beta
+{1\over 2}\,(\bar\phi\alpha+\bar\psi\beta)(\phi\bar\alpha-\psi\bar\beta)\\
&+(\bar\phi\alpha+\bar\psi\beta)\biggl\{
(\phi\bar\alpha-\psi\bar\beta)(1-x_1x_2)+\phi\bar\beta x x_2-\psi\bar\alpha\bar x x_1\biggr\}\,.\\
\end{aligned}$$

Now taking the conjugate of $-\begin{pmatrix}\alpha\\\beta\end{pmatrix}$ and 
multiplying on the left with $\begin{psmallmatrix}0&-1\\1&0\\\end{psmallmatrix}$, we find
$$\begin{aligned}
\begin{pmatrix}\bar\beta\\-\bar\alpha\end{pmatrix}&=
\begin{pmatrix} x_1&\bar x\\x&x_2\\\end{pmatrix}\\
&\biggl [(1+\phi\psi)\begin{pmatrix}-\bar\psi\\\bar\phi\\\end{pmatrix}
+{1\over 2}(1+\bar\phi\bar\psi)(\phi\bar\alpha+\psi\bar\beta)
\begin{pmatrix}\alpha\\\beta\end{pmatrix}\biggr],
\end{aligned}$$
and so
$$\begin{aligned}
-\begin{pmatrix}\alpha\\\beta\\\end{pmatrix}&=
(1+\bar\phi\bar\psi)\begin{pmatrix}x_2&-\bar x\\-x&x_1\\\end{pmatrix}\begin{pmatrix} \phi\\\psi\end{pmatrix}\\
&+{1\over 2}(1+\phi\psi)(\bar\phi\alpha+\bar\psi\beta)
\begin{pmatrix}x_2&-\bar x\\-x&x_1\\\end{pmatrix}\begin{pmatrix}x_1&\bar x\\x&x_2\\\end{pmatrix} \\
&\hskip .5in\biggl [
(1+\phi\psi)\begin{pmatrix}-\bar\psi\\\bar\phi\\\end{pmatrix}+{1\over 2}(1+\bar\phi\bar
\psi)(\phi\bar\alpha+\psi\bar\beta)\begin{pmatrix}\alpha\\\beta\\\end{pmatrix}
\biggr ]\\
&=(1+\bar\phi\bar\psi)\begin{pmatrix}x_2&-\bar x\\-x&x_1\\\end{pmatrix}\begin{pmatrix} \phi\\\psi\end{pmatrix}\\
&+{1\over 2}(1+\phi\psi)(\bar\phi\alpha+\bar\psi\beta)
(x_1x_1-x\bar x)\\
&\hskip .5in\biggl [
(1+\phi\psi)\begin{pmatrix}-\bar\psi\\\bar\phi\\\end{pmatrix}+{1\over 2}(\phi\bar\alpha+\psi\bar\beta)\begin{pmatrix}\alpha\\\beta\\\end{pmatrix}
\biggr ]\\
&=(1+\bar\phi\bar\psi)\begin{pmatrix}x_2&-\bar x\\-x&x_1\\\end{pmatrix}\begin{pmatrix} \phi\\\psi\end{pmatrix}\\
&+{1\over 2}\biggl[ \bar\phi\bar\psi(1+\phi\psi)\begin{pmatrix}\alpha\\\beta\\\end{pmatrix}+\alpha\beta(\phi\bar\alpha+\psi\bar\beta)\begin{pmatrix}-\bar\psi\\\bar\phi\\\end{pmatrix}\biggr]\\
&=(1+\bar\phi\bar\psi)\begin{pmatrix}x_2&-\bar x\\-x&x_1\\\end{pmatrix}\begin{pmatrix} \phi\\\psi\end{pmatrix}+{1\over 2} \bar\phi\bar\psi(1+\phi\psi)\begin{pmatrix}\alpha\\\beta\\\end{pmatrix},
\end{aligned}$$
where the last equality relies on the earlier expression for the product $\alpha\beta$.  It follows that
$$\biggl[1+{1\over 2}\bar\phi\bar\psi(1+\phi\psi)
\biggr]\begin{pmatrix}\alpha\\\beta\end{pmatrix}=-(1+\bar\phi\bar\psi)\begin{pmatrix}x_2&-\bar x\\-x&x_1\\\end{pmatrix}\begin{pmatrix} \phi\\\psi\end{pmatrix},$$
and so
$$\begin{pmatrix}\alpha\\\beta\\\end{pmatrix}=(1+{1\over 2}\bar\phi\bar\psi)
\begin{pmatrix}-x_2&\bar x\\x&-x_1\\\end{pmatrix}\begin{pmatrix}\phi\\\psi\\\end{pmatrix},$$
as was asserted.

In fact, we find that $\alpha\beta=\phi\psi$, and a further rewarding computation
involving our solution for $\alpha,\beta$ and $\vartheta$ confirms that the
coordinates $x_1,x_2,x$ are indeed all scaled by the factor
$K=1+{1\over 2}(\phi\psi+\bar\phi\bar\psi)+{3\over 4}\phi\psi\bar\phi\bar\psi$
under the transformation $u(\alpha,\beta)$ for the specified $\alpha$ and $\beta$.\end{proof}

In the notation of the lemma, $$u.(x_1,x_2,x|\phi,\psi)=K(\phi,\psi)\,(x_1,x_2,x|0,0),$$ and in fact one finds
$K^{-2}=(1-\phi\psi-\bar\phi\bar\psi)$.


\begin{theorem}\label{thm:vol}
For $K=K(\phi,\psi)=1+{1\over 2}(\phi\psi+\bar\phi\bar\psi)+{3\over 4}\phi\psi\bar\phi\bar\psi$, the three-form
$$\begin{aligned}
{\rm Vol}&={\rm dlog}_e\, (K^{-1} x_2)\wedge {\rm d}(K^{-1}v)\wedge{\rm d} (K^{-1} u)\\
\end{aligned}$$
in $(x_1,x_2,u+iv{\rm |}\phi,\psi)$ coordinates on $\hyperp\subset{\mathbb R }^{3,1|4}$
is invariant under the Hermitean action of ${\rm OSp}_{\mathbb C}(1{\rm |} 2)$, and the body 
of ${\rm Vol}$ is the hyperbolic volume form.  Moreover, ${\rm Vol}={\rm d}\Omega$ admits the primitive
$$\Omega=
{{(1-\phi\psi-\bar\phi\bar\psi)}\over 2}\biggl\{{\rm d}v\wedge{\rm d}u+{\rm d log}_{e}x_2\wedge[v{\rm d}u-u{\rm d}v]\biggr\}\,.$$
\end{theorem}

\begin{proof}
As in \cite{defect}, there is
an equivariant tower of isometric actions
$$\begin{aligned}
&{\rm OSp}_{\mathbb C}(1|2)&\acts&\hskip 5ex\hyperp\\
&\hskip 3ex\vee&&\hskip 5ex\cup\\
&\hskip 2.5ex\widehat{\rm SL}(2,{\mathbb C})&\acts&\hskip 5ex\widehat{\hyperp}\\
&\hskip 3ex\vee&&\hskip 5ex\cup\\
&{\rm SL}(2,{\mathbb C})&\acts&\hskip 5ex\widecheck\hyperp\, ,\\
\end{aligned}$$
where the bottom of the tower is Wigner's action on the classical
hyperboloid.  The classical projection from $\widecheck\hyperp$
to the upper half-space model is given by
$(x_1,x_2,u+iv)\mapsto {2\over {x_2}}(u,v,1)$.
The invariance of the classical hyperbolic volume
form
$$-{\rm d log}_e\, x_2\wedge {\rm d}u\wedge{\rm d}v=-{i\over 2}\,\,{\rm d}{\rm log}_e\, x_2\wedge{\rm d} x\wedge{\rm d}\bar x, ~{\rm for}~x=u+iv,$$
under the action ${\rm SL}(2,{\mathbb C})\acts\widecheck\hyperp$ tautologoically
implies invariance on $\widehat{\hyperp}$
under $\widehat{\rm SL}(2,{\mathbb C})$.

The previous lemma provides a mapping
$\hyperp\to\widehat{\hyperp}$ from the hyperboloid in ${\mathbb R}^{3,1|4}$
to its purely bosonic part given by $$(x_1,x_2,x|\phi,\psi)\mapsto K(\phi,\psi)(x_1,x_2,x|0,0),$$
where $K(\phi,\psi)=1+{1\over 2}(\phi\psi+\bar\phi\bar\psi)+{3\over 4}\phi\psi\bar\phi\bar\psi$,
and pulling back the $\widehat{\rm SL}(2,{\mathbb C})$-invariant form under this mapping immediately gives the asserted formula for ${\rm Vol}$.

Continuing to compute, we have
$$\begin{aligned}
{\rm d}(K^{-1}w)&=K^{-2}(K\,{\rm d} w-w\,{\rm d}K), {\rm ~for}~w\in\{u,v\},\\ 
{\rm dlog}_e(K^{-1}x_2)&={\rm dlog}_e x_2-{\rm dlog}_e K,\\
\end{aligned}$$
and it follows that
$$
{\rm Vol}=K^{-2}\,{\rm dlog}_e\,x_2\wedge {\rm d}v\wedge{\rm d} u+X\wedge{\rm d}Y$$
where
$$\begin{aligned}
X&={1\over 2}\,K^{-2},~{\rm and}\\
Y&={\rm d}v\wedge{\rm d} u+{\rm dlog}_e\,x_2\wedge(  v\,{\rm d}u-u\,{\rm d} v).\\
\end{aligned}$$
Thus for any finite-sided polyhedron ${\mathbb T}$, 
$$\begin{aligned}
\int_{\mathbb T} {\rm Vol}&=K^{-2}\int_{\mathbb T} {\rm dlog}_e x_2\,{\rm d}v\,{\rm d}u+\int_{\mathbb T}X{\rm d}Y\\
&=K^{-2}\int_{\mathbb T} {\rm dlog}_e x_2\,{\rm d}v\,{\rm d}u
-\int_{\mathbb T}Y{\rm d}X+\int_{\partial{\mathbb T}}XY\\
&=\int_{\partial{\mathbb T}}XY,\\
\end{aligned}$$
where the second equality relies on Stokes' Theorem and the third on a fortuitous cancellation. \end{proof}

\section{Triangles and Tetrahedra in Super Space}\label{tris}

\subsection{Triangles}

Let $\triangle$ denote a super triangle in $\hyperp$ with non-collinear vertices
$$\begin{aligned}
{P}&=(p_1,p_2,p|\alpha,\beta),\\
{Q}&=(q_1,q_2,q|\gamma,\delta),\\
{R}&=(r_1,r_2,r|\varepsilon,\varphi),\\
\end{aligned}$$
where the subscripted coordinates are real and the others complex. As follows from Corollary \ref{cor:tgt}, $\triangle$ satisfies the Hyperbolic Law of Cosines.

Define the  
three even variables $$d=\langle P,Q\rangle, ~e=\langle Q,R\rangle,~f=\langle R,P\rangle$$
and consider the function
$$H(d,e,f)=2def+1-d^2-e^2-f^2.$$
In fact, $H(d,e,f)>0$, as follows from the Hyperbolic Law of Cosines\footnote{Heron's formula gives the tangent of half the hyperbolic area of $\triangle$
as ${{\sqrt{H(d,e,f)}}\over{1+d+e+f}}$, but this formula does not hold for the area of a super triangle
due to the fermionic correction in Theorem 6.2 of \cite{defect}, cf. Section \ref{sec:close}. This nevertheless proves
$H>0$.}.

\begin{lemma}\label{real}
Given any triple $P,Q,R\in{\mathbb R}^{3,1|4}$ of points in the notation above so that no
pairwise difference is isotropic, we can arrange by
applying an element of ${\rm OSp}_{\mathbb C}(1|2)$
that  $p,q,r\in{\mathbb C}$ have a common real part as well as the existence
of fermions $\mu,\rho,\sigma,\tau$ and $A\in{\rm SL}(2,{\mathbb C})$
so that 
$\begin{psmallmatrix}\alpha\\\beta\\\end{psmallmatrix}=A\begin{psmallmatrix}\mu\\\rho\\\end{psmallmatrix}$,
$\begin{psmallmatrix}\gamma\\\delta\\\end{psmallmatrix}=A\begin{psmallmatrix}\mu\\\sigma\\\end{psmallmatrix}$,
and $\begin{psmallmatrix}\varepsilon\\\varphi\\\end{psmallmatrix}=A\begin{psmallmatrix}\mu\\\tau\\\end{psmallmatrix}$.
\end{lemma}

\noindent Notice that there remains one further degree of freedom among the bosonic coordinates
in this parametrization of orbits of triples.
Moreover,
analogues of Lemma \ref{real} likewise hold for other classical normalizations, for example for
super triangles in $\hyperp$.  Notice that in dimension two, there is an additional constraint on orbits of triples
coming from their ordered bodies giving a positive orientation.

\begin{proof} 
First perturb $\triangle$ with an element of $\widehat{\rm SL}(2,{\mathbb C})$ so that
$$\begin{aligned}
t&=(\bar qp_1-\bar pq_1)+(\bar rq_1-\bar qr_1)+(\bar pr_1-\bar rp_1)\\
&=(p_1-q_1)(\bar q-\bar r)-(q_1-r_1)(\bar p-\bar q)\\
&=p_1(\bar q-\bar r)+q_1(\bar r-\bar p)+r_1(\bar p-\bar q)\\
\end{aligned}$$
has non-zero body.
Given a pair of  fermions $\xi,\eta$, let us suppose that $u(\xi,\eta). P=(p_1',p_2',p'|\alpha',\beta'),$
and likewise respectively $\gamma'$ and $\varepsilon '$ for ${Q}$ and ${R}$, where
$u(\xi,\eta)\in{\rm OSp}_{\mathbb C}(1|2)$ as in Lemma \ref{factor}.

According to the second identity in Corollary \ref{corgood}, we have
$$\begin{aligned}
\alpha'&=
p_1\xi+\bar p\eta+(1+\xi\eta+{1\over 2}\xi\eta\bar\xi\bar\eta)\alpha,\\
\gamma'&=
q_1\xi+\bar q\eta+(1+\xi\eta+{1\over 2}\xi\eta\bar\xi\bar\eta)\gamma,\\
\varepsilon'&=r_1\xi+\bar r\eta+(1+\xi\eta+{1\over 2}\xi\eta\bar\xi\bar\eta)\varepsilon,\\
\end{aligned}$$
so $\alpha'=\gamma'=\varepsilon'$ is equivalent to the linear system
$$\begin{pmatrix}
p_1-q_1&\bar p-\bar q\\
q_1-r_1&\bar q-\bar r\\
\end{pmatrix}
\begin{pmatrix} \xi\\\eta\\\end{pmatrix}=\begin{pmatrix}\gamma-\alpha\\\varepsilon-\gamma\\\end{pmatrix}$$
with determinant $t\neq 0$, as required.  We have thus proved

\begin{lemma}\label{ace}
Given any generic triple of points in ${\mathbb R}^{3,1|4}$, there are odd $\xi,\eta$
and $u=u(\xi,\eta)\in{\rm OSp}_{\mathbb C}(1|2)$ as in Lemma \ref{factor} so that the
image under $u$ of the triple share a common first fermionic coordinate.\hfill\qed
\end{lemma}

\noindent We think of one fermion (the second coordinate) associated to each point,
together with another odd ``Manin type invariant'' (the common first coordinate) associated to 
the triple itself.

Thus without loss of generality, we may hence assume that $\alpha=\gamma=\varepsilon$ and shall apply
$g=\begin{psmallmatrix} a&b&0\\c&d&0\\0&0&1\\\end{psmallmatrix}\in\widehat{SL}(2,{\mathbb C})$.
According to the third identity in Corollary \ref{corgood}, the bosonic coordinates of $g.P$ are given
by
$$\begin{pmatrix}p_1'&\bar p'\\p&p_2'\\\end{pmatrix}
=\begin{pmatrix}\bar a&\bar c\\\bar b&\bar d\end{pmatrix}
\begin{pmatrix}p_1&\bar p\\p&p_2\\\end{pmatrix}
\begin{pmatrix}a&b\\c&d\\\end{pmatrix}$$
and likewise for $Q$ and $R$.  

It follows that for $x,y\in\{ p,q,r\}$, we have
$$\begin{pmatrix}x_1'-y_1'&\bar x'-\bar y'\\x'-y'&x'_2-y'_2\\\end{pmatrix}
=\begin{pmatrix}\bar a&\bar c\\\bar b&\bar d\end{pmatrix}
H(x,y)
\begin{pmatrix}a&b\\c&d\\\end{pmatrix},$$
where the Hermitean matrix
$$H(x,y)=\begin{pmatrix}x_1-y_1&\bar x-\bar y\\x-y&x_2-y_2\\\end{pmatrix}$$
has
determinant 
$\langle X-Y,X-Y\rangle$ using that $X,Y$
have a common first fermionic coordinate, and this determinant is non-zero by hypothesis.

Thus, $x$ and $y$ have a common real part if and only if
$$0=\begin{pmatrix}\bar b&\bar d\\\end{pmatrix} H(x,y)\begin{pmatrix}a\\c\\\end{pmatrix}
+\begin{pmatrix}\bar a&\bar c\\\end{pmatrix} H(x,y)\begin{pmatrix}b\\d\\\end{pmatrix}.$$
There is a simultaneous solution 
$g=\begin{psmallmatrix} a&b&0\\c&d&0\\0&0&1\\\end{psmallmatrix}\in\widehat{SL}(2,{\mathbb C})$ 
for pairs $(x,y)=(p,q)$ and $(x,y)=(q,r)$
exactly as in the bosonic case (namely, write $a,b,c,d\in{\mathbb C}$ in polar coordinates and diagonalize one
of the Hermitean matrices over the unitary group).

Again by the second equation in Corollary \ref{corgood}, we find that the fermions $\begin{psmallmatrix} \alpha\\\beta\\\end{psmallmatrix}$, $\begin{psmallmatrix} \gamma\\\delta\\\end{psmallmatrix}$, and
$\begin{psmallmatrix} \varepsilon\\\varphi\\\end{psmallmatrix}$, with $\alpha=\gamma=\varepsilon$
each transform under $g\in{\rm OSp}_{\mathbb C}(1|2)$ by
$A=\begin{psmallmatrix}\bar a&\bar c\\\bar b&\bar d\\\end{psmallmatrix}$, giving the final assertion.
\end{proof}

Turning now to the ideal case, $P,Q,R\in{\mathcal L}_+$ are the  {\it ideal vertices} of an {\it ideal super triangle} $\triangle$ which we interpret as follows.  According to Theorem \ref{thm:lines}, a super geodesic with a distinguished point gives rise to a pair
$E,F\in{\mathcal L}_+$ with $\langle E,F\rangle =2$, where $U={1\over 2} (E+F)$ is the distinguished point.  Altering $E,F$ in their rays (that is, scaling by a boson with positive body) preserving this condition
corresponds to translating the distinguished point along the geodesic.  Regarding $\triangle
=\triangle({P,Q,R})\subseteq\hyperp$ as a triple of geodesics pairwise sharing asymptotic rays in ${\mathcal L}_+$, $\triangle$ is interpreted as a triple of distinct rays in ${\mathcal L}_+$, the rays of $P,Q,R\in{\mathcal L}_+$
or points in them the ideal vertices.

\begin{lemma} \label{satis}
Given any three pairwise distinct rays in ${\mathcal L}_+$, there are unique respective points $X_i$ in these rays so that $\langle X_i,X_j\rangle=2$, for distinct $i$ and $j$ in $\{1,2,3\}$.
\end{lemma}

\begin{proof} Let $Y_i\in{\mathcal L}_+$ denote any respective points in the three rays, for $i=1,2,3$, and define $X_i=\sqrt{2\,{\langle Y_j,Y_k\rangle}\over{\langle Y_i,Y_j\rangle \langle Y_i,Y_k\rangle}}$, for $\{i,j,k\}=\{1,2,3\}$.
\end{proof}

\begin{lemma} \label{parm}
For any $P,Q,R\in{\mathcal L}_+$ satisfying the conclusions of the previous lemma, the ideal triangle
$\triangle=\triangle({P,Q,R})$ is parametrized by
$$X(s,t)=2{P \over t}+{t\over s}\,\biggl\{{Q\over{s_-}}+{R\over{s_+}}-2{P\over {s}}
\biggr\},~{\rm where}~s_\pm=s\pm\sqrt{s^2-4},$$
for $2<s<\infty$ and $0<t<{1\over 2}{{s_+}}$.
\end{lemma}

\begin{proof} According to Theorem \ref{thm:lines}, the super geodesic asymptotic to the rays of $Q$ and $R$ is parametrized
by $X=X(s)={1\over 2}\, (sQ+s^{-1}R )$, so we find $\langle P,X\rangle =s+s^{-1}=\hat s$.  
Now using Corollary \ref{useful} and writing $Y=tP+(1-t)X$, we compute
$$\langle Y, Y\rangle =(1-t)^2+2t(1-t)\hat s=t^2(1-2\hat s)+2t(\hat s-1)+1,$$
so $\langle Y, Y\rangle =0$ for $t=(1-2\hat s)^{-1}$, whence $(1-2\hat s)\,Y=P-\hat s(sQ+s^{-1}R)$. 

Now computing $Z=cY$ with $2=\langle P,Z\rangle$, we find that $c={{2\hat s-1}\over \hat s^2}$, whence
$Z=-{P\over \hat s^2}+{{s^2Q+R}\over {s^2+1}}$, and it follows from a second application of
Theorem \ref{thm:lines} that 
$$2X(s,t)=\biggl({{\hat s^2-t^2}\over{t\hat s^2}} \biggr ) P
+\biggl({{t}\over{1+s^{-2}}} \biggr ) Q
+\biggl({{t}\over{1+s^{2}}} \biggr ) R,
$$
for $0<s<\infty$, $0<t<\hat s$.
Since $\hat s=s+s^{-1}=2\,{\rm cosh}\,{\rm log}_es$, we find
$s^{\pm 1}={1\over 2}\biggl(\hat s\pm\sqrt{\hat s^2-4}\biggr )$
using ${\rm cosh}^{-1} d={\rm log}_e(d+\sqrt{d^2-1})$, and hence
$1+s^{\pm 2}={{\hat s}\over 2}\biggl (\hat s\pm\sqrt{\hat s^2-4}
\biggr ),$
so the result follows upon replacing the parameter $\hat s$ with $s$.\end{proof}


\subsection{Tetrahedra}


Four points $P_1,P_2,P_3,P_4\in\hyperp^3$, no three of which are collinear, determine a (possibly degenerate) {\it super tetrahedron}
${\mathbb T}\subseteq \hyperp^3$.  Define the inner products 
$$d_{ij}=\langle P_i,P_j\rangle={\rm cosh}\, D(P_j,P_j),~{\rm for}~i,j=1,2,3,4,$$
where $D$ denotes the distance function,
and organize them into the {\it Gram edge matrix} $G({\mathbb T})_{ij}=(d_{ij}),$ for $i,j=1,2,3,4$.
The usual formula \cite{junand,junand2} for dihedral angles in terms of the Gram matrix is unchanged
for super tetrahedra:

\begin{lemma}\label{dih}
The cosine $a_{ij}$ of the dihedral angle of ${\mathbb T}$ along the edge opposite $\overline{P_iP_j}$
    is given by $a_{ij}=-{{c_{ij}}\over {\sqrt{c_{ii}}\sqrt{c_{jj}}}},$
    where the $c_{ij}$ are the cofactors of the Gram edge matrix, where $i,j=1,2,3,4$.
\end{lemma}

\begin{proof} We shall prove the formula explicitly for the edge $\overline{P_3P_4}$, and to this end,
pick some $X\in\overline{P_3P_4}$.  According to Corollary \ref{cor:tgt}, the unit vector
$v$ parallel to $\overrightarrow{XP_3}$ is given by
$$v={{P_3-X\langle P_3,X\rangle}\over{\sqrt{\langle P_3,X\rangle^2-1}}},$$
and likewise the respective unit vectors parallel to $\overrightarrow{XP_1}$ and $\overrightarrow{XP_2}$
are
$$u_i={{P_i-X\langle P_i,X\rangle}\over{\sqrt{\langle P_i,X\rangle^2-1}}},~{\rm for}~i=1,2.$$
Thus, $w_i=u_i-\langle u_i,v\rangle\, v$ is perpendicular to $\overline{P_3P_4}$ of length
$1-\langle u_i,v\rangle^2$.

Taking now $X=P_4$, we find $\langle u_i,v\rangle=d_{i3}-d_{i4}d_{34}$, for $i=1,2$, and
it follows from elementary algebra that the cosine of dihedral angle is
$$\begin{aligned}
&{{\langle w_1,w_2\rangle}\over{\sqrt{\langle w_1,w_1\rangle}
\sqrt{\langle w_2,w_2\rangle}}}
={{\langle u_1,u_2\rangle-\langle u_1,v\rangle\,\langle u_2,v\rangle}\over{
\sqrt{1-\langle u_2,v\rangle^2}
\sqrt{1-\langle u_1,v\rangle^2}}}\\
&\quad\quad\quad={{(d_{34}^2-1)(d_{12}-d_{14}d_{24})+
(d_{23}-d_{24}d_{34})
-(d_{13}-d_{14}d_{34})}\over
{\sqrt{H(d_{13},d_{14},d_{34})}
\sqrt{H(d_{23},d_{24},d_{34})},
}}\\
\end{aligned}$$
which we finally recognize in terms of cofactors of the Gram matrix.
The assertion for an arbitrary edge of ${\mathbb T}$ follows by symmetry.
\end{proof}

The ${\rm OSp}_{\mathbb C}(1|2)$ moduli space of tetrahedra evidently has real dimension
$4(3|4)-2(3|2)=(6|12)$ with moduli given by the 6 inner products of vertex pairs plus the
four-dimensional span of the Dirac fermions on a face, together with the two-dimensions
of Dirac fermions on its opposite vertex.  
On the other hand, applying Lemma \ref{real} we may assume
that a tetrahedron has vertices 
$$\begin{aligned}
P&=(p_1,p_2,u^*+iv^p|(\mu~\rho)\,A^t)\\
Q&=(q_1,q_2,u^*+iv^q|(\mu~\sigma)\,A^t),\\
R&=(r_1,r_2,u^*+iv^r|(\mu~\tau)\,A^t),\\
W&=(w_1,w_2,u^w+iv^w|\xi,\eta),\\
\end{aligned}$$
with real bosons, complex fermions and $A\in{\rm SL}(2,{\mathbb C})$.

As shown in \cite{defect}, the area of a typical ideal super triangle diverges,
and we analogously 
have

\begin{theorem}\label{diverges}
Consider the ideal super tetrahedron $\mathbb T$ with vertices $P,Q,R,W\in{\mathcal L}_+$ 
in the notation above.  Then the volume of $\mathbb T$ diverges provided at least one of
$\mu\rho$, $\mu\sigma$, or $\mu\tau$ is non-zero.
\end{theorem}

\begin{proof} 
Suppose it is $\mu\rho\neq 0$. 
As before for ideal triangles, we may assume that $P,Q,R$ satisfy the conclusions of
Lemma \ref{satis} and hence adopt the parametrization $X(s,t)$  given in Lemma \ref{parm} for the ideal triangle
$\triangle=\triangle(P,Q,R)$.

By Stokes' Thorem and Theorem \ref{thm:vol}, the volume of ${\mathbb T}$ is given by
$$\int_{\mathbb T} {\rm Vol}=\int_{\partial {\mathbb T}}\Omega=\int_{\partial {\mathbb T}}
{{(1-\phi\psi-\bar\phi\bar\psi)}\over 2}\biggl\{{\rm d}v\wedge{\rm d}u+{\rm d log}_{e}x_2\wedge[v{\rm d}u-u{\rm d}v]\biggr\},$$
and of course $\int_\triangle \Omega$ is one of four summands in this expression, corresponding to the four codimension-one faces of $\mathbb T$.  We shall prove by direct computation that $\int_\triangle \Omega$ diverges.
To this end since the $u$-coordinate is constant on $\Delta$ by construction, we find
$$\int_{\partial {\mathbb T}}\Omega={{u^*}\over 2}\,\int_{\partial {\mathbb T}}(1-\phi\psi-\bar\phi\bar\psi)\,
{\rm d log}_{e}x_2\wedge{\rm d}v.$$

According to Lemma \ref{parm},
$$
x_2(s,t)=2{{p_2}\over t}+{t\over s}\,X_2(s),~{\rm where}~X_2(s)=
{{q_2}\over{s_-}}+{{r_2}\over{s_+}}-2\,{{p_2}\over s},$$
so
$${\rm d}x_2=\biggl (t^2X_2-2s\,p_2\biggr ){{{\rm d}t}\over {st^2}}
+\biggl (sX_2'-X_2\biggr ){{t\,{\rm d}s}\over{s^2}},
$$
where $X_2'={{{\rm d}X_2}\over{{\rm d} s}}$.
Likewise for 
$$v(s,t)=2{{v^p}\over t}+{t\over s}\,V(s),~{\rm where}~V(s)=
{{v^q}\over{s_-}}+{{v^r}\over{s_+}}-2\,{{v^p}\over s},$$
so that
${\rm d log}_{e}x_2\wedge{\rm d}v={{{\rm d}x_2}\over {x_2}}\wedge{\rm d}v$ is given by
$$
\biggl [{{(t^2X_2-2s\,p_2)(sV'-V)-(t^2V-2s\,v^p)(sX_2'-X_2)}
\over{t^2X_2+2s\,p_2}}\biggr ]\, {{{\rm d}t\wedge{\rm d}s}\over {s^2}}.
$$

Similarly  according to Lemma \ref{parm}, we find
$$\begin{pmatrix} \phi(s,t)&\psi(s,t)\\\end{pmatrix}A^{-t}={{s^2-t^2}\over {s^2t}}\,
\begin{pmatrix} \mu&\rho\\\end{pmatrix}+t\,\begin{pmatrix} \mu&
2\{{{\sigma}\over{s+\sqrt{s^2-4}}}
+
{{\tau}\over{s-\sqrt{s^2-4}}}\}
\\\end{pmatrix}.
$$
Since 
$\begin{pmatrix}\xi'&\eta'\\\end{pmatrix}=\begin{pmatrix} \xi&\eta\\\end{pmatrix}B$
implies that $\xi\eta=\xi'\eta'$ for $B\in{\rm SL}(2,{\mathbb C})$,
we conclude that
$$\phi\psi=(t+{1\over t}-{t\over{s^2}})\,\mu\,\biggl[
\biggl({{s^2-t^2}\over {s^2t}}\biggr)\rho+2t\biggl(
{\sigma\over{s+\sqrt{s^2-4}}}+{\tau\over{s-\sqrt{s^2-4}}}
\biggr)
\biggr].$$
It is only the coefficient ${1\over{t^2}}+1-{2\over{s^2}}+{{t^2}\over{s^4}}(1-s^2)$ of $\mu\rho$ in $\phi\psi$
that concerns us, and indeed only its lowest-order summand ${1\over {t^2}}$ in $t$.

Thus, the coefficient of $\mu\rho$ in $\int_\triangle \Omega$ contains a summand 
$$\int_2^\infty \int_0^{{1\over 2}(s+\sqrt{s^2-4})} {1\over t^2}\,
\biggl[{{t^2 C+D}\over{t^2X_2+2s\,p_2}}\biggr]
 {{{\rm d}t\,{\rm d}s}\over {s^2}},$$
where $C,D,X_2$ are typically non-vanishing functions of $s$ independent of $t$.
It therefore follows that this coefficient diverges.
\end{proof}

\section{Closing Remarks}\label{sec:close}

The classical Angle Defect Theorem follows from the Schl\"afli formula.
See Danny Calegari's lovely blog \cite{danny} for proofs and a contextual discussion
of the classical Schl\"afli formula and \cite{luo} for relevant computational details.
Yet the Angle Defect Theorem fails \cite{defect} for super triangles, so the
general Schl\"afli formula as it stands must be false for super tetrahedra.
It is an interesting and challenging open problem to formulate and compute
the fermionic correction to the Schl\"afli formula, and there are hints:
  Milnor's proof 
in \cite{danny} together with Lemma \ref{real}, and
the general scheme of \cite{defect}  using Lemma \ref{dih} are both suggestive in this regard.
It is also interesting to wonder what aspects of \cite{wpt} might extend to the 
super case in light of Theorem \ref{diverges}.

The difference between the area and the angle
defect provides a novel invariant of super triangles which is
additive under disjoint union.  In the same spirit and going
beyond Hilbert's Third Problem in three-dimensional hyperbolic space, we wonder what are the
additive invariants of hyperbolic super tetrahedra.

As mentioned in the Introduction, \cite{others} posits three metrics on
${\mathcal N}=1$ super upper half-space which are unrelated to
Hermitean conjugation and our considerations here.  The authors 
probe several interesting questions including: Are the metrics Einstein?
What are their Green's function? What is the connection with Arakelov geometry
\cite{yuri}?  We might equally well ask all these same questions for
our metric here. 

We mentioned 
in Section \ref{sec:vol} that the classical projection from $\widecheck\hyperp^3$
to the upper half-space model ${\mathcal U}$ is given by
$(x_1,x_2,u+iv)\mapsto {2\over {x_2}}(u,v,1)$.  
Following
\cite{ahlfors}
and introducing the imaginary quaternionic unit $j$, 
the nice fact is that the map $(x_1,x_2,u+iv)\mapsto q={2\over{x_2}}(u+iv+j)$
from $\check\hyperp^3$\ to ${\mathcal U}$
taking the height in half-space
as the coordinate of $j$
 is equivariant for the action of $g=\begin{psmallmatrix}
a&b\\c&d\end{psmallmatrix}\in{\rm SL}(2,{\mathbb C})$ by $A\mapsto g^\dagger A g$ on Minkowski space and by fractional linear transformation $q\mapsto{{aq+b}\over{cq+d}}$ on this quaternionic version of $\mathcal U$.
Moreover, \cite{others} argues that one takes not only $q\in{\mathcal U}$
but also a quaternionic $j\phi+\psi$ and then acts by
super fractional linear transformations.
In any case, we ask what is the natural equivariant map
from $\hyperp^3$ to complex ${\mathcal N}=1$ super upper half-space? 

It is the quadratic equation 
$$\begin{aligned}
0&=\vartheta^2[\alpha\beta+\bar\alpha\bar\beta+2\alpha\beta\bar\alpha\bar\beta]+\vartheta[X+(Y-\bar Y)+2(Z-\bar Z)]\\
&+X(Y-\bar Y)+{1\over 2}(Y-\bar Y)^2+{1\over 2}X^2-2\alpha\beta\bar\alpha\bar\beta(\phi\psi+\bar\phi\bar\psi)
\end{aligned}$$ 
in $\vartheta$, where $
X=x_1\bar\alpha\alpha+x_2\bar\beta\beta+x\bar\beta\alpha +\bar x\bar\alpha\beta$, $
Y=\bar\alpha\phi+\bar\beta\psi$,
and $Z=\alpha\beta Y$,
that implies invariance of our metric.
Meanwhile, quadratic equations with Grassmann coefficients defy usual 
bosonic expectations, for instance with no quadratic formula for roots if the discriminant
or the leading coefficient have vanishing body, both of which happen here.  Such equations can have any number of solutions. 
It follows that our extension of Wigner's representation based upon the solution
$\vartheta=-{1\over 2}[X+(Y-\bar Y)+(Z-\bar Z)]$ may not be unique,
though we can say that it is the unique solution in the linear span of $X,Y,Z$.
There is also the assumption in certain lemmas in the Appendix that $\vartheta$ is pure imaginary, without
which the algebra becomes onerous and whose relaxation 
provides another quadratic equation and potentially different extension.
As the quadratic in $\vartheta$ above figures so prominently, Igor Frenkel asks if the considerations of this paper might be related to the affine algebras  \cite{igor}
of type $B(O,\ell)$ since the latter also determine a quadratic equation
 ``nearly'' determined from a root system.

Finally, the matrices
we consider here
are not quite Hermitean since the diagonal entry $\vartheta$ is pure imaginary rather
than real.  This suggests various generalizations of Hermitean matrices
that may be worth considering, specifically in the super case.  It also suggests that the algebraically
messy alternative mentioned in the previous paragraph that $\vartheta$  is not pure imaginary may warrant
further investigation.

\appendix

\section{Proof of the Main Theorem}

\begin{theorem}\label{thm} 
Setting $A=\begin{psmallmatrix}x_1&\bar x&\phi\\x&x_2&\psi\\
-\bar\phi& -\bar\psi&\vartheta&\end{psmallmatrix}$, the action 
$A\mapsto g^\dagger A g$ with $g^\dagger$ the~conjugate~super~transpose~of $g\in{\rm OSp}_{\mathbb C}(1|2),$ preserves the
quadratic form $Q(x_1,x_2,x{\rm |}\phi,\psi)=x_1x_2-x\bar x+\phi\psi+\bar\phi\bar\psi$,
where we have defined $\vartheta=-{1\over 2}[X+(Y-\bar Y)+(Z-\bar Z)]$,
with $
X=x_1\bar\alpha\alpha+x_2\bar\beta\beta+x\bar\beta\alpha +\bar x\bar\alpha\beta$, $
Y=\bar\alpha\phi+\bar\beta\psi$,
and $Z=\alpha\beta Y$.

\end{theorem}

\begin{proof}
One computes directly that 
$g^\dagger A g
=\begin{psmallmatrix}x_1'&\bar{x}'&\phi'\\
x'&x_2'&\psi'\\
-\bar{\phi}&-\bar{\psi}&\vartheta'\end{psmallmatrix}$ 
with
\medskip
\begin{equation*}
\begin{aligned}[t]
{x_1'}&=A+A_1+A_2, ~{\rm where}\\
A&=a\bar{a}x_1+c\bar{c}x_2+(a\bar{c}x+\bar{a}c\bar{x}),\\
A_1&=(a\bar{\gamma}\bar{\phi}+\bar{a}\gamma\phi)+(c\bar{\gamma}\bar{\psi}+\bar{c}\gamma\psi),\\
A_2&=-\vartheta\bar{\gamma}\gamma;\\
\end{aligned}
\qquad
\begin{aligned}[t]
{x_2'}&=B+B_1+B_2, ~{\rm where}\\
B&=b\bar{b}x_1+d\bar{d}x_2+(b\bar{d}x+\bar{b}d\bar{x}),\\
B_1&=(b\bar{\delta}\bar{\phi}+\bar{b}\delta\phi)+(d\bar{\delta}\bar{\psi}+\bar{d}\delta\psi),\\
B_2&=-\vartheta\bar{\delta}\delta;\\
\end{aligned}
\end{equation*}
\begin{equation*}
\begin{aligned}[t]
x'&=C+C_1+C_2, {\rm where}\\
C&=a\bar{b}x_1+c\bar{d}x_2+a\bar{d}x+\bar{b}c\bar{x},\\
C_1&=\bar\delta(a\bar\phi+c\bar\psi)+\gamma(\bar b\phi+\bar d\psi);\\
C_2&=-\vartheta\bar{\delta}\gamma;\\
\end{aligned}
\qquad
\begin{aligned}[t]
\vartheta'&=W_1\alpha+W_2\beta+fW,{\rm where}\\
W_1&=x_1\bar{\alpha}+x\bar{\beta}+\bar{f}\bar{\phi},\\
W_2&=x_2\bar{\beta}+\bar{x}\bar{\alpha}+\bar{f}\bar{\psi},\\
W&=\bar{\alpha}\phi+\bar{\beta}\psi+\vartheta \bar{f};\\
\end{aligned}
\end{equation*}

\begin{equation*}
\begin{aligned}[t]
\phi'&=U_1\alpha+U_2\beta+fU, {\rm where}\\
U_1&=\bar{a}x_1+\bar{c}x+\bar{\gamma}\bar{\phi},\\
U_2&=\bar{c}x_2+\bar{a}\bar{x}+\bar{\gamma}\bar{\psi},\\
U&=\bar{a}\phi+\bar{c}\psi+\vartheta\bar{\gamma};\\
\end{aligned}
\qquad
\begin{aligned}[t]
\psi'&=V_1\alpha+V_2\beta+fV, {\rm where}\\
V_1&=\bar{b}x_1+\bar{d}x+\bar{\delta}\bar{\phi},\\
V_2&=\bar{d}x_2+\bar{b}\bar{x}+\bar{\delta}\bar{\psi},\\
V&=\bar{b}\phi+\bar{d}\psi+\vartheta\bar{\delta}.\\
\end{aligned}
\end{equation*}

It follows that 
$$\begin{aligned}
x_1'x_2'-x'\bar{x}'&=[AB-C\bar{C}]+[AB_1+BA_1-C\bar{C}_1-\bar{C}C_1]\\
&+[A_1B_1-C_1\bar{C}_1]+[A_2B_2-C_2\bar{C}_2]\\
&+[A_2(B+B_1)+B_2(A+A_1)-C_2(\bar{C}+\bar{C}_1)-\bar{C}_2(C+C_1)].\\
\end{aligned}$$
Routine computations using the defining relations
for ${\rm OSp}_{\mathbb C}(1|2)$ 
imply

\begin{lemma} \label{lemq}We have the following identities
$$\begin{aligned}
AB-C\bar{C}&=(x_1x_2-x\bar{x})(1-\alpha\beta)(1-\bar{\alpha}\bar{\beta}),\\
AB_1+BA_1-C\bar{C}_1-\bar{C}C_1&=x_1\{\alpha(\bar{\alpha}\bar{\beta}-1)\psi+\bar{\alpha}(\alpha\beta-1)\bar{\psi}\}\\
&+x_2\{\beta(1-\bar{\alpha}\bar{\beta})\phi+\bar{\beta}(1-\alpha\beta)\bar{\phi}\}\\
&+x\{\alpha(1-\bar{\alpha}\bar{\beta})\phi+\bar{\beta}(\alpha\beta-1)\bar{\psi}\}\\
&+\bar{x}\{\bar{\alpha}(1-{\alpha}{\beta})\bar{\phi}+{\beta}(\bar{\alpha}\bar{\beta}-1){\psi}\},\\
A_1B_1-C_1\bar{C}_1&=2(1-\alpha\beta)\bar{\alpha}\bar{\phi}\bar{\beta}\bar{\psi}+2(1-\bar{\alpha}\bar{\beta})\alpha\phi\beta\psi\\
&-\{\alpha\bar{\alpha}\bar{\phi}\phi+\beta\bar{\beta}\bar{\psi}\psi+\beta\bar{\alpha}\bar{\psi}\phi+\alpha\bar{\beta}\bar{\phi}\psi\},\\
A_2B_2-C_2\bar{C}_2&=\alpha\bar{\alpha}\beta\bar{\beta}\vartheta(\vartheta-\bar\vartheta),\\
A_2(B+B_1)+B_2(A+A_1)&-C_2(\bar{C}+\bar{C}_1)-\bar{C}_2(C+C_1)\\
&=-\vartheta\{\bar{\alpha}\alpha x_1+\bar{\beta}\beta x_2+\bar{\beta}\alpha x+\bar{\alpha}\beta\bar{x}\}\\
&-2\vartheta\{\alpha\beta(\bar{\alpha}\phi+\bar{\beta}\psi)-\bar{\alpha}\bar{\beta}(\alpha\bar{\phi}+\beta\bar{\psi})\}\\
&+\delta\bar{\gamma}(C+C_1)(\vartheta+\bar{\vartheta}).\\
\end{aligned}$$
\end{lemma}
\noindent Notice that the body of the identity for $AB-C\bar{C}$ is due to Wigner.

It likewise follows that
$$\begin{aligned}
\phi'\psi'&=f^2UV+\alpha\beta(U_1V_2-U_2V_1)+f\{(UV_1-VU_1)\alpha+(UV_2-VU_2)\beta\}\\
&=f^2UV+\alpha\beta(U_1V_2-U_2V_1)+\{(UV_1-VU_1)\alpha+(UV_2-VU_2)\beta\}\\
\end{aligned}$$
since $f=1+\alpha\beta$.
Again routine computations imply

\begin{lemma} \label{lemp} We have the {following} identities
$$\begin{aligned}
UV_1-VU_1&=(x\phi-x_1\psi)(1-\bar{\alpha}\bar{\beta})+\bar{\phi}(\bar{\alpha}\phi+\bar{\beta}\psi)
+\vartheta\{ x_1\bar{\alpha}+x\bar{\beta}+2\bar\alpha\bar\beta\bar{\phi}\},\\
UV_2-VU_2&=(x_2\phi-\bar{x}\psi)(1-\bar{\alpha}\bar{\beta})+\bar{\psi}(\bar{\alpha}\phi+\bar{\beta}\psi)+\vartheta\{x_2\bar{\beta}+\bar{x}\bar{\alpha}+2\bar\alpha \bar\beta\bar{\psi}\},\\
U_1V_2-U_2V_1&=(x_1x_2-x\bar{x})(1-\bar{\alpha}\bar{\beta})-2\bar{\gamma}\bar{\delta}\bar{\phi}\bar\psi
+(x_2\bar{\beta}+\bar{x}\bar{\alpha})\bar{\phi}-(x_1\bar\alpha+x\bar\beta)\bar \psi,\\
UV&=(1-\bar{\alpha}\bar{\beta})\phi\psi+\vartheta(\bar\alpha \phi+\bar\beta \psi)+\vartheta^2\bar{\alpha}\bar{\beta}.\\
\end{aligned}$$
\end{lemma}

Given $g=\begin{psmallmatrix}a&b&\alpha\\c&d&\beta\\\gamma&\delta&f\\\end{psmallmatrix}\in{\rm OSp}_{\mathbb C}(1{\rm |}2)$
and $(x_1,x_2,x|\phi,\psi)\in{\mathbb R}^{3,1|4}$,
define the quantities
$$\begin{aligned}
X&=x_1\bar\alpha\alpha+x_2\bar\beta\beta+x\bar\beta\alpha +\bar x\bar\alpha\beta, \\
Y&=\bar\alpha\phi+\bar\beta\psi,\\
Z&=\alpha\beta Y.\\
\end{aligned}$$
Notice that each of $X, (Y-\bar Y), (Z-\bar Z)$ is pure imaginary and that
$\alpha\beta X=0=\bar\alpha\bar\beta X$, $\alpha\beta\bar Y=0=\bar\alpha\bar\beta Y$, and
$X^2=-2\alpha\beta\bar\alpha\bar\beta(x_1x_2-x\bar x)$.



Let us now combine these calculations and turn to the proof that
the action 
$\hat A\mapsto g^\dagger \hat A g$ preserves the
quadratic form $Q(x_1,x_2,x{\rm |}\phi,\psi)=x_1x_2-x\bar x+\phi\psi+\bar\phi\bar\psi$
for $\vartheta=-{1\over 2}[X+(Y-\bar Y)+(Z-\bar Z)]$.
Indeed, provided $\bar\vartheta=-\vartheta$ we find that
$\phi\psi+\bar\phi\bar\psi$ is given by

$$\begin{aligned}
&\boxed{\bigl [{\color{red} (1-\bar\alpha\bar\beta)(x\phi-x_1\psi)}+{\color{orange}\bar\phi(\bar\alpha\phi+\bar\beta\psi)}+\vartheta(x_1\bar\alpha+x\bar\beta+2\bar\alpha\bar\beta\bar\phi\bigr)]\alpha}\\
+&\bigl [{\color{red}(1-\alpha\beta)(\bar x\bar \phi-x_1\bar \psi)}+{\color{orange}\phi(\alpha\bar\phi+\beta\bar\psi)}-\vartheta(x_1\alpha+\bar x\beta+2\alpha\beta\phi\bigr)]\bar\alpha\\
+&\boxed{\bigl [{\color{red}(1-\bar\alpha\bar\beta)(x_2\phi-\bar x\psi)}+{\color{orange}\bar\psi(\bar\alpha\phi+\bar\beta\psi)}+\vartheta(x_2\bar\beta+\bar x\bar\alpha+2\bar\alpha\bar\beta\bar\psi\bigr)]\beta}\\
+&\bigl [{\color{red} (1-\alpha\beta)(x_2\bar\phi- x\bar\psi)}+{\color{orange}\psi(\alpha\bar\phi+\beta\bar\psi)}-\vartheta(x_2\beta+ x\alpha+2\alpha\beta\psi\bigr)]\bar\beta\\
+&\boxed{\bigl[{\color{blue}(x_1x_2-x\bar x)(1-\bar\alpha\bar\beta)-2\bar\alpha\bar\beta\bar\phi\bar\psi}+{\color{green}(x_2\bar\beta+\bar x\bar\alpha)\bar\phi
-(x_1\bar\alpha+x\bar\beta)\bar\psi}\bigr ]\alpha\beta}\\
+&\bigl[{\color{blue}(x_1x_2-x\bar x)(1-\alpha\beta)-2\alpha\beta\phi\psi}+{\color{green}(x_2\beta+ x\alpha)\phi
-(x_1\alpha+\bar x\beta)\psi}\bigr]\bar\alpha\bar\beta\\
+&\boxed{(1+2\alpha\beta)\bigl[{\color{blue}(1-\bar\alpha\bar\beta)\phi\psi}+\vartheta(\bar\alpha\phi+\bar\beta\psi)+\vartheta^2\bar\alpha\bar\beta\bigr]}\\
+&(1+2\bar\alpha\bar\beta)\bigl[{\color{blue}(1-\alpha\beta)\bar\phi\bar\psi}-\vartheta(\alpha\bar\phi+\beta\bar\psi)+\vartheta^2\alpha\beta\bigr],\\
\end{aligned}$$
where the boxed terms express $\phi'\psi'$ directly from 
Lemma \ref{lemp}, and the unboxed terms likewise express $\bar\phi'\bar\psi'$.

Direct elementary computation provides the following identities.

\begin{lemma}\label{lem14} Assume that $\bar\vartheta=-\vartheta$.  Then in the expression above for $\phi'\psi'+\bar\phi'\bar\psi'$, we have
$$\begin{aligned}
{\rm the~sum~of~red~terms}&=-(AB_1+BA_1-C\bar{C}_1-\bar{C}C_1),\\
{\rm the~sum~of~orange~terms}&=-2\,Y\bar Y,\\
{\rm the~sum~of~green~terms}&=X(Y-\bar Y),\\
\hskip 1.5em{\rm the~sum~of~black~terms}&=\vartheta^2[(1+2\alpha\beta)\bar\alpha\bar\beta+(1+2\bar\alpha\bar\beta)\alpha\beta]\\
&+\vartheta[2X+(Y-\bar Y)+4(Z-\bar Z)].\\
\end{aligned}$$
\end{lemma}
\noindent The cancellation engendered by this double appearance of the red terms
in $\phi'\psi'+\bar\phi'\bar\psi'$ and in the second identity of Lemma \ref{lemq} provides the motivating 
insight for the computations here.

Now, $x_1'x_2'-x'\bar x'$ is the sum of the expressions in Lemma \ref{lemq},
and so upon applying Lemma \ref{lem14}, the grand total $x_1'x_2'-x'\bar x'+\phi'\psi'+\bar\phi'\bar\psi'$ is given by
$$\begin{aligned}
&(x_1x_2-x\bar x)[1-\alpha\beta+\alpha\beta\bar\alpha\bar\beta+{\color{blue}
(1-\bar\alpha\bar\beta)\alpha\beta+(1-\alpha\beta)\bar\alpha\bar\beta}]\\
&-2[(1-\alpha\beta)\bar\alpha\bar\beta\bar\phi\bar\psi
+(1-\bar\alpha\bar\beta)\alpha\beta\phi\psi+{\color{blue}\alpha\beta\bar\alpha\bar\beta(\phi\psi+\bar\psi\bar\psi)}]\\
&+{\color{blue}(1+2\alpha\beta)((1-\bar\alpha\bar\beta)\phi\psi+
+(1+2\bar\alpha\bar\beta)(1-\alpha\beta)\bar\phi\bar\psi}\\
&+Y\bar Y+{\color{blue} X(Y-\bar Y)-2Y\bar Y}\\
&+\vartheta\biggl [
-X-2(Z-\bar Z)+{\color{blue}2X+(Y-\bar Y)
+4(Z-\bar Z)}
\biggr ]\\
&+\vartheta^2\biggl [
{\color{blue}\alpha\beta(1+2\bar\alpha\bar\beta)+\bar\alpha\bar\beta(1+2\alpha\beta)}-2\alpha\beta\bar\alpha\bar\beta
\biggr],\\
\end{aligned}$$
where the blue terms arise from (all colors of terms of) $\phi'\psi'+\bar\phi'\bar\psi'$ from above, and the black ones arise from
$x_1'x_2'-x'\bar x'$.  Simplifying the difference
$$\Delta=x_1'x_2'-x'\bar x'+\phi'\psi'+\bar\phi'\bar\psi'-
(x_1x_2-x\bar x+\phi\psi+\bar\phi\bar\psi),$$ 
we find the quadratic in $\vartheta$ given by

$$\begin{aligned}
\Delta&=\vartheta^2[\alpha\beta+\bar\alpha\bar\beta+2\alpha\beta\bar\alpha\bar\beta]+\vartheta[X+(Y-\bar Y)+2(Z-\bar Z)]\\
&+X(Y-\bar Y)-Y\bar Y-\alpha\beta\bar\alpha\bar\beta(x_1x_2-x\bar x)\\
&-(\bar\alpha\bar\beta+ 
2\alpha\beta\bar\alpha\bar\beta)\phi\psi
-(\alpha\beta+ 
2\alpha\beta\bar\alpha\bar\beta)\bar\phi\bar\psi\\
&=\vartheta^2[\alpha\beta+\bar\alpha\bar\beta+2\alpha\beta\bar\alpha\bar\beta]+\vartheta[X+(Y-\bar Y)+2(Z-\bar Z)]\\
&+X(Y-\bar Y)+{1\over 2}(Y-\bar Y)^2+{1\over 2}X^2-2\alpha\beta\bar\alpha\bar\beta(\phi\psi+\bar\phi\bar\psi),
\end{aligned}$$
using $X^2=-2\alpha\beta\bar\alpha\bar\beta(x_1x_2-x\bar x)$ and
$Y^2=2\bar\alpha\phi\bar\beta\psi$ in the last line.

In this form, one might guess that a solution may be a linear combination of $X,(Y-\bar Y),(Z-\bar Z)$
and then confirm that this specific linear combination
$$\vartheta=-{1\over 2}[X+(Y-\bar Y)+(Z-\bar Z)],$$
solves the equation $\Delta=0$, as required, through a remarkable cancellation or two.
(In order to simplify the summand in $\Delta$ which is quadratic in $\vartheta$, recall
 that multiplication by $\alpha\beta$ or $\bar\alpha\bar\beta$ annihilates
 $X$ or $Y\bar Y$.)  
Directly comparing with the formula for $\vartheta'$ at the beginning of this section, one finds that furthermore $\vartheta'=-2\vartheta$, which is of unclear signficance. This completes the proof of Theorem \ref{thm} and the Main Theorem from the Introduction.
\end{proof}

\begin{remark}The function 
$$\begin{aligned} 
({\mathbb R}^{3,1|4})^{\times 2}&\to\hat {\mathbb C}[0]\\
A\times B&\mapsto\kappa(A,B)=-2{{\rm d}\vartheta(A,B)}\\
\end{aligned}$$
is given by
$$\kappa(\hat A_1,\hat A_2)=
(\bar\psi_2\,\phi_1-\psi_2\,\bar\phi_1)+(\phi_2\,\bar\psi_1-\bar\phi_2\,\psi_1)$$
in the obvious notation.  $\kappa$ is a skew-symmetric  pairing on ${\mathbb R}^{3,1|4}$ taking purely imaginary values and satisfying 
$\kappa(iA, B)=-\kappa(A,iB)$.  
Indeed,
taking the exterior derivative of $\vartheta$ and evaluating at the origin, one 
finds that 
$$\begin{aligned}
-2{\rm d}\vartheta(\hat A_1,\hat A_2)&={\rm d}[(\bar\alpha\phi+\bar\beta\psi)-(\alpha\bar\phi+\beta\bar \psi)]\\
&=
{\rm d}\bar\alpha\,\phi-{\rm d}\alpha\,\bar\phi+{\rm d}\bar\beta\,\psi-{\rm d}\beta\,\bar\psi.\\
\end{aligned}$$Direct evaluation
immediately gives the asserted formula for $\kappa$, and the further assertions about it
are readily verified.  
Presumably $\kappa$ is
the skew part of the Hermitean form on the complex Lie super algebra
arising from the Killing form on the compact real form as in the Introduction.
\end{remark}

\clearpage
\setcounter{section}{0}
\setcounter{equation}{0}
\setcounter{theorem}{0}

\begin{center}
{\LARGE\bfseries Corrigendum to {\it 3d super hyperbolic geometry}\par}
\vspace{4pt}
{\large with an Appendix correcting {\it Angle defect for super triangles}\par}
\vspace{12pt}
{\large Robert Penner\par}
\vspace{4pt}
Institut des Hautes \'Etudes Scientifiques\\
35 route de Chartres, 91440 Bures-sur-Yvette, France\\
\texttt{rpenner@ihes.fr}
\end{center}
\vspace{12pt}
\begin{center}
\small
\textbf{Corrigendum for the \emph{Journal of Geometry and Physics}}\\
to R.\ C.\ Penner, \emph{3d super hyperbolic geometry},
J.\ Geom.\ Phys.\ \textbf{199} (2024), 105171.
\end{center}

\begin{quote}
\small
\noindent\textbf{Abstract.} \noindent We correct two errors in the volume calculation of
\emph{3d super hyperbolic geometry}, together with an incidental formula
for the transformed auxiliary matrix entry $\vartheta'$ in its Appendix.
The Main Theorem, the quadratic form, the invariant metric, and the
super-hyperboloid construction are correct and unchanged.  What fails is
the identification of a bosonizing ordinary three-form with intrinsic
super-volume: the printed form contains $K^{-1}$ in place of $K$, its
printed primitive has the wrong sign, and the resulting three-form is not
invariant under the odd action.  Following the metric--Berezinian viewpoint
for hyperbolic superspaces emphasized by Hu and Zong, a direct computation
from the invariant $3|4$ metric gives a finite intrinsic Berezinian volume
for an ideal super tetrahedron.  The published cusp-divergence calculation
remains correct for the non-invariant three-form actually used there, but
not as a statement about intrinsic volume.  The appendix records the
parallel correction in two dimensions.
\end{quote}

\section*{Introduction}

The principal purpose of this corrigendum is to correct the volume discussion
in Sections~4--5 of \cite{corrPenner3d}, published in the \emph{Journal of
Geometry and Physics}.  The Main Theorem, giving the modified
Hermitian-conjugation action of $\OSp_{\mathbb C}(1|2)$ and the invariance of
its quadratic form, is correct, as are the resulting super hyperboloid and
its invariant metric.  The error occurs when that metric geometry is
replaced by an ordinary differential three-form obtained from the
bosonizing normalization.  In Theorem~4.2, the normalization factor is also
written as $K^{-1}$ where the pullback calculation gives $K$, and the
printed primitive has the wrong overall sign.  More fundamentally, even the
corrected bosonizing pullback is not equivariant for the odd action and
therefore is not the intrinsic invariant volume of the supermanifold.

The appropriate intrinsic object is instead the Berezinian density determined
by the invariant supermetric.  This is the general metric construction used
for hyperbolic superspaces by Hu and Zong \cite{corrHuZong}; their model has
different odd dimension, so no formula from their paper is imported here.
Rather,  the Berezinian density is computed directly from the $3|4$ metric of \cite{corrPenner3d},
demonstrating that its top odd coefficient is a constant multiple of the ordinary
hyperbolic volume density.  Consequently an ideal tetrahedron has finite
intrinsic Berezinian volume.  The local divergence calculation of Section~5
nevertheless remains algebraically correct for the particular non-invariant
three-form used there; what fails is its interpretation as super-hyperbolic
volume.

The same distinction corrects the earlier two-dimensional paper
\cite{corrDefect}.  The appendix here records that, for a triangle already in
a suitable gauge, the calculations in that paper remain exactly correct
for the printed two-form, though they are unrelated to intrinsic metric
area.  Curiously, that same form is closely related to a differential form
arising independently in the subsequent theory of additive supertriangle
invariants \cite{corrGP}.

\section{The Main Theorem in three dimensions}

The Main Theorem of \cite{corrPenner3d} is correct as stated, and the
calculation in the Appendix proving invariance of the quadratic form $
        Q(x_1,x_2,x\mid\phi,\psi)
        =
        x_1x_2-x\bar x+\phi\psi+\bar\phi\,\bar\psi
$ under the modified Hermitian conjugation is correct. 
Specifically with
the auxiliary parameter $\vartheta$ defined there, the transformation
$
        \widehat A\longmapsto g^\dagger \widehat A g$, for 
   $g\in\OSp_{\mathbb C}(1|2)
$
induces an action on $\mathbb R^{3,1|4}$ preserving $Q$,
and consequently, the super hyperboloid defined by $Q=1$, its metric induced by this
quadratic form, and the results of Sections~1--3 of \cite{corrPenner3d} are
unchanged.

There is, however, one incidental error at the end of the Appendix of
\cite{corrPenner3d} which is not used elsewhere in the paper, as follows.
After the proof of invariance of $Q$ is complete, it is inaccurately
stated that direct comparison with the previously displayed transformation
formula gives
\[
        \vartheta'=-2\vartheta,
\]
where
in the notation of \cite{corrPenner3d}, 
with
$
 g=\begin{psmallmatrix}
 a&b&\alpha\\ c&d&\beta\\ \gamma&\delta&f
 \end{psmallmatrix}
 \in \OSp_{\mathbb C}(1|2)
$
and
$$
 X=x_1\bar\alpha\alpha+x_2\bar\beta\beta
   +x\bar\beta\alpha+\bar x\bar\alpha\beta,
\qquad
 Y=\bar\alpha\phi+\bar\beta\psi,
$$
the auxiliary parameter is
\begin{equation}\label{eq:vartheta-def-corr}
 \vartheta
 =-\frac{f\bar f}{2}\bigl[X+(Y-\bar Y)\bigr].
\end{equation}
On the other hand, the explicit matrix multiplication in the Appendix gives
\begin{equation}\label{eq:vartheta-prime-first}
 \vartheta'
 =X+fY-\bar f\,\bar Y+f\bar f\,\vartheta.
\end{equation}
Since $f=1+\alpha\beta$, \eqref{eq:vartheta-def-corr} and \eqref{eq:vartheta-prime-first} actually give
\begin{equation}\label{eq:vartheta-prime-corr}
 {\vartheta'=(f\bar f-2)\vartheta}=-\frac{\vartheta}{f\bar f}
\end{equation}
where the second equality follows from $\alpha\beta\bar\alpha\bar\beta\,\vartheta=0$.

This correction has no effect on the Main Theorem.  Nor does it affect the subsequent skew pairing
\[
 \kappa=-2\,d\vartheta
\]
defined in the final remark of the Appendix.  The factor $-2$ in that
definition is a normalization of $d\vartheta$ and is unrelated to the erroneous statement about $\vartheta'$.
Furthermore, the pairing \(\kappa\) is indeed, not just {\sl presumably} as stated in \cite{corrPenner3d}  up to the conventional normalization and sign of the Killing form, the ordinary skew-symmetrization of the odd-odd part of the super-Hermitian form on \(\mathfrak{osp}(1|2;\mathbb C)\) induced by the Killing form of its compact real form \(UOSp(1|2)\).

\section{Correction of the three-dimensional volume calculation}

There are two concrete errors in Theorem~4.2 and its proof.  Lemma~4.1
correctly sends
\[
 (x_1,x_2,x\mid\phi,\psi)
 \longmapsto
 K(\phi,\psi)(x_1,x_2,x\mid0,0),
\]
where
\[
 K
 =1+\frac12(\phi\psi+\bar\phi\,\bar\psi)
 +\frac34\phi\psi\bar\phi\,\bar\psi,
 \qquad
 K^{-2}=1-\phi\psi-\bar\phi\,\bar\psi.
\]
Thus, the literal pullback of classical hyperbolic volume under this
bosonizing normalization is
\begin{equation}\label{eq:3d-correct-vol}
 \Vol_{\rm pb}
 =d\log(Kx_2)\wedge d(Kv)\wedge d(Ku),
\end{equation}
not the expression printed in Theorem~4.2  with $K^{-1}$ replacing $K$ in \eqref{eq:3d-correct-vol}.   Independently,
the displayed primitive there has the wrong overall sign.  A correctly
signed primitive of \eqref{eq:3d-correct-vol} is
\begin{equation}\label{eq:3d-correct-primitive}
\Omega_{\rm pb}
=-\frac12\left\{
d(Kv)\wedge d(Ku)
+d\log(Kx_2)\wedge
\bigl[Kv\,d(Ku)-Ku\,d(Kv)\bigr]
\right\},
\end{equation}
so that $d\Omega_{\rm pb}=\Vol_{\rm pb}$.  The more important point is that
\eqref{eq:3d-correct-vol} is only a bosonizing pullback.  The normalization
is not equivariant for the odd action, so this ordinary three-form is not the
intrinsic invariant super-volume.

As emphasized in the general theory of hyperbolic superspaces by Hu and Zong
\cite{corrHuZong}, an invariant supermetric canonically determines a Berezinian
volume density.  Their explicit models are not used here; for the present
$3|4$ space, the calculation is made directly from the metric of
\cite{corrPenner3d}.  Write
\[
 x=U+iV,\qquad
 \phi=\alpha+i\beta,\qquad
 \psi=\gamma+i\delta,
 \qquad
 \epsilon=\alpha\gamma-\beta\delta,
\]
so
\[
 Q=x_1x_2-U^2-V^2+2\epsilon,
 \qquad
 x_2=\frac{1+U^2+V^2-2\epsilon}{x_1}.
\]
A direct Berezinian calculation gives
\begin{equation}\label{eq:3d-berezinian}
{
 \sqrt{\operatorname{Ber}g}
 =\frac1{x_1}\bigl(1-2\epsilon-2\epsilon^2\bigr).}
\end{equation}
The degree-zero coefficient is the ordinary hyperbolic volume density.  Since
$\epsilon^2=2\alpha\beta\gamma\delta$, the top odd coefficient is
$-4/x_1$ for the displayed odd ordering, hence a constant multiple of the
same classical volume density.  In particular the invariant metric
Berezinian introduces no additional cusp singularity.

\section{Ideal super tetrahedra}

The geometric parametrization and local asymptotic calculations in
Section~5 of \cite{corrPenner3d} remain algebraically correct for the ordinary
three-form used there.  In particular, the displayed $t^{-2}$ singularity at
an ideal vertex is genuinely present for that form.  What changes is its
geometric interpretation.  Because the three-form is obtained from the
non-equivariant bosonizing construction rather than from the invariant
metric, that divergence is not a divergence of intrinsic super-hyperbolic
volume.

By the intrinsic super-volume of a superdomain we mean the Berezin integral of the Berezinian density determined by the invariant supermetric.
For this intrinsic Berezinian volume, the conclusion is the opposite.
By
\eqref{eq:3d-berezinian}, Berezin integration over the four odd directions
leaves, up to the orientation convention, a fixed multiple of ordinary
hyperbolic volume.  Since an ordinary ideal hyperbolic tetrahedron has finite
volume, so does the corresponding ideal super tetrahedron with respect to
the intrinsic Berezinian density.  The infinite-volume conclusion of
Theorem~5.6 thus does not survive as a statement about metric volume.

This is exactly parallel to the two-dimensional correction below: the
published cusp calculation remains a correct calculation for the particular
ordinary form under consideration, while the intrinsic area or volume is obtained by 
Berezin integration of the density determined directly by the invariant metric.

\appendix

\section{Correction to {\it Angle defect for super triangles}}
\label{app:2d}

There are two points to correct in \cite{corrDefect}.  First, the sign
$1-\phi\psi$ printed in Lemma~4.2 is a typographical error: the bosonic
normalization is
$
B(x_1,x_2,y\mid\phi,\psi)
=(1+\phi\psi)(x_1,x_2,y\mid0,0),
$
as is already used correctly later in the paper.  Second and
substantively, the two-form
\begin{equation}\label{eq:oldform}
\Omega_{\rm old}
=
\frac{
d[(1+\phi\psi)x_2]\wedge d[(1+\phi\psi)x_1]
}{2(1+\phi\psi)y}
\end{equation}
is not invariant under the full odd action of $\osp$.  It is the pullback
of ordinary hyperbolic area under the above bosonizing normalization, but
that normalization is not $\osp$-equivariant.

For a triangle already normalized as in Lemma~6.1,
all of the subsequent computations in Sections~5--7 are exactly correct:
the remaining normalizations are even, so no erroneous use of odd
$\osp$-invariance enters.  What changes is their geometric interpretation.
They are computations for the particular form $\Omega_{\rm old}$, not for
an intrinsic super-area form.  In particular, the explicit degree-two
discrepancy computed in \cite{corrDefect} is genuinely present and is nonzero in
general.  In the notation of that paper,
\[
\mathcal D_{\rm sup}-\mathcal A_{\rm old}
=
\frac{\mathcal P+\mathcal Q}{\sqrt{1-\mathcal R^2}},
\]
so the super Angle Defect Theorem remains false.  Likewise, the divergence of
$\int\Omega_{\rm old}$ for a generic normalized ideal triangle is a correct calculation, although it is unrelated to the
intrinsic area question.

There is also a curious and potentially significant coincidence.  The form
$\Omega_{\rm old}$, originally introduced here in an attempted area
interpretation, is closely related to a differential form that arises
independently in the new theory of regular additive invariants of spin
supertriangles \cite{corrGP}.  Thus, the form and normalized calculations
built from it appear to encode genuine supergeometric information even
though $\Omega_{\rm old}$ is not the invariant area form.  We do not pursue
that relation here, other than to mention that $\mathcal D_{\rm sup}-\mathcal A_{\rm old}$ 
is in fact an invariant of supertriangles which is additive under subdivision into two supertriangles.

By the intrinsic super-area of a superdomain we mean the Berezin integral of the 
Berezinian density determined by the invariant supermetric.
For this intrinsic Berezinian area, the invariant metric provides the appropriate density.
For
\[
Q=x_1x_2-y^2+2\phi\psi=1
\]
with the metric induced by the invariant bilinear form of \cite{corrDefect}, its
Berezinian density in local coordinates $(s,y\mid\phi,\psi)$, $s=x_1$, is
\begin{equation}\label{eq:2d-berezinian}
\sqrt{\operatorname{Ber}g}
=
\frac1s\bigl(1-2\phi\psi\bigr),
\end{equation}
up to the sign determined by the ordering of the odd coordinates.  The
degree-zero coefficient is thus the ordinary hyperbolic area density, while the
top odd coefficient is its constant multiple $\pm2$.  Both are integrable
over an ideal hyperbolic triangle: the former gives $\pi$ and the latter
$\pm2\pi$.  Hence the invariant metric Berezinian has no cusp divergence.

Accordingly, the divergence found in \cite{corrDefect} should not be interpreted
as divergence of intrinsic super-area.  What survives from the original
argument is the normalized calculation for $\Omega_{\rm old}$ and its
nonzero degree-two correction; what replaces the claimed invariant area form
is the finite Berezinian density \eqref{eq:2d-berezinian} determined directly
by the invariant metric.

\end{document}